\documentclass[preprint,12pt]{elsarticle}

\usepackage[pdftex]{hyperref}
\biboptions{sort&compress,numbers,square}

\usepackage{algorithmic}
\usepackage{amsfonts}
\usepackage{graphicx}
\usepackage{epsfig}
\usepackage{color}
\usepackage{geometry}
\usepackage{tikz}
\usetikzlibrary{trees,positioning,shapes}
\usepackage{subcaption} 
\usepackage{amsmath}
\usepackage{amsthm}

\usepackage{xcolor, colortbl}

\usepackage{float}
\usepackage[hypcap]{caption}
\usepackage{subcaption}
\usepackage{algorithm2e}
\usepackage{array}

\newcommand{\revision}[1]{\textcolor{black}{#1}}
\newcommand{\inblue}[1]{\textcolor{black}{#1}}
\newcommand{\inred}[1]{\textcolor{black}{#1}}

\newcommand{\ingreen}[1]{\textcolor{black}{#1}}

\usepackage{bm}
\usepackage{mathtools}

\DeclareMathOperator*{\argmin}{arg\,min}

\usepackage{rotating}

\begin{document}

\begin{frontmatter}

\title{Automatic stabilization of finite-element simulations using neural networks and hierarchical matrices}

\author{Tomasz S\l{}u\.zalec$^{(1)}$, Mateusz Dobija$^{(1)}$, Anna Paszy\'nska$^{(1)}$, Ignacio Muga$^{(2)}$,  \\ Maciej Paszy\'{n}ski$^{(3)}$}

\address{$^{(1)}$Jagiellonian University, Krak\'ow, Poland \\ $^{(2)}$ Pontificia Universidad Cat\'olica of Valpara\'iso, Chile \\  $^{(3)}$ AGH University of Science and Technology,
Krak\'{o}w, Poland \\
e-mail: maciej.paszynski@agh.edu.pl }

\begin{abstract}
Petrov-Galerkin formulations with optimal test functions allow for the stabilization of finite element simulations. In particular, given a discrete trial space, the optimal test space induces a numerical scheme delivering the best approximation in terms of a problem-dependent energy norm.
This ideal approach has two shortcomings: first, we need to explicitly know the set of optimal test functions; and second, the optimal test functions may have large supports inducing expensive dense linear systems. 
 
Nevertheless, parametric families of PDEs are an example where it is worth investing some (offline) computational effort to obtain stabilized linear systems that can be solved efficiently, for a given set of parameters, in an online stage.
Therefore, as a remedy for the first shortcoming, we explicitly compute (offline) a function mapping any PDE-parameter, to the matrix of coefficients  of optimal test functions (in a basis expansion) associated with that PDE-parameter.
Next, as a remedy for the second shortcoming, we use the low-rank approximation to hierarchically compress the (non-square) matrix of coefficients of optimal test functions. In order to accelerate this process, we train a neural network to learn a critical bottleneck of the compression algorithm (for a given set of PDE-parameters).
When solving \emph{online} the resulting (compressed) Petrov-Galerkin formulation, we employ a GMRES iterative solver with inexpensive matrix-vector multiplications thanks to the low-rank features of the compressed matrix. We perform experiments showing that the full online procedure \inred{as fast as} the original (unstable) Galerkin approach. 
We illustrate our findings by means of 2D Eriksson-Johnson and Hemholtz model problems.
\end{abstract}
	

\end{frontmatter}

\tableofcontents

\section{Introduction}
Difficult finite-element simulations solved with the Galerkin method 
(where we employ the same trial and test space) often generate incorrect numerical results with oscillations or spurious behavior. Examples of such problems are the advection-dominated diffusion equation~\cite{Erikkson} or the Helmholtz equation~\cite{Helmholtz}.

Petrov-Galerkin formulations\footnote{i.e., trial and test spaces are not equal, although they share the same dimension.}~\cite{PG} with optimal test functions allow for automatic stabilization of difficult finite-element simulations. 
This particular Petrov-Galerkin approach is equivalent to the residual minimization (RM) method~\cite{Erikkson}, whose applications include advection-diffusion \cite{iGRM1,iGRM2}, Navier-Stokes \cite{iGRM3}, or space-time formulations \cite{ResMinSpaceTime}.

In general, for a fixed trial space, RM allows for stabilization by enriching the discrete test space where optimal test functions of the associated Petrov-Galerkin method live 
\footnote{Contrary to standard stabilization methods like SUPG~\cite{SUPG,SUPG2}, there are no special stabilization terms modifying the weak formulation.}. 
%
If we had explicitly the formulas for the optimal test functions expanded in the enriched discrete test space, then we can return to the Petrov-Galerkin formulation and solve problems in the best possible way for a given trial space.

This scenario has two shortcomings. The first problem is that the computation of optimal test functions is expensive. It requires solving a large system of linear equations\footnote{As large as the dimension of the enriched discrete test space.} with multiple right-hand sides (one right-hand side per each basis function of the trial space).
The second problem is that the optimal test functions can have global supports, and thus the Petrov-Galerkin method with optimal test functions can generate a dense matrix, expensive to solve. 

Nevertheless, parametric families of partial differential equations (PDEs) are an example where it is worth investing some (offline) computational effort, to obtain stabilized linear systems that
can be solved efficiently, for a given set of parameters in an online stage.
Therefore, in the context of a parametric family of PDEs, and as a remedy for
the first shortcoming, we explicitly compute (offline) a function mapping any PDE-parameter, to the matrix of coefficients of optimal test functions\footnote{expanded in the basis of the enriched discrete test space.} associated with that PDE-parameter. We emphasize that this last procedure
is independent of particular right-hand side sources and/or prescribed boundary data of the PDE family. 

The obtained matrix $\mathbb W$ of optimal test functions coefficients is dense.
The Petrov-Galerkin method induces a linear system of the form $\mathbb B^T\mathbb W\,x=(L^T\mathbb W)^T$, where $L$ is a right-hand side vector and $\mathbb B$ is the matrix associated with the bilinear form of the underlying PDE.
To avoid dense matrix computations, we compress $\mathbb W$ using the approach of hierarchical matrices \cite{Hackbush}.
Having the matrix $\mathbb W$ compressed into a hierarchical matrix $\mathbb H$, we employ the GMRES method \cite{Saad}, which involves computations of the residual $R = \mathbb B^T\mathbb H\,x-(L^T\mathbb H)^T$
and the hierarchical matrix enables matrix-vector multiplication of $\mathbb H\,x$ and $L^T\mathbb H$ in a quasi-linear computational cost.

However, compressing the matrix $\mathbb W$ for each PDE-parameter is an expensive procedure. Thus, with the help of an artificial neural network, we train (offline) the critical bottleneck of the compression algorithm. 
\inred{We obtain a stabilized method with the additional quasi-linear cost resulting from matrix-vector multiplications within the GMRES method.
From our numerical results with the Ericksson-Johnson and Hemholtz model problems, this cost of stabilization (the cost of compression of the hierarchical matrix and the cost of GMRES with  hierarchical matrix multiplication by a vector) is of the same order as the cost of the solution of the original Galerkin problem without the stabilization. Thus, we claim we can obtain stabilization with neural networks and hierarchical matrices practically for free.}

\subsection{One-dimensional illustration of stabilization}
Let us illustrate how we stabilize difficult computational problems by means of a one-dimensional example for the advection-dominated diffusion model. 

Given $0<\epsilon\leq 1$, consider the following differential equation:
\begin{equation}\label{eq:AdvReac_1d}
\left\{
\begin{array}{rl}
-\epsilon u''+u'=\inred{0} & \hbox{in } (0,1)\,; \\
-\epsilon u'(0)+u(0)=1 & \hbox{and } u(1)=0.  
\end{array}\right.
\end{equation}
%
In weak-form, equation~\eqref{eq:AdvReac_1d} translates into find $u\in H_{0)}^1(0,1):=\{v\in H^1(0,1): v(1)=0\}$ such that:
\begin{equation}\label{eq:varfor}
\epsilon\int_0^1 u'v' + \int_0^1 u'v +u(0)v(0) = \inred{v(0)} 
\end{equation}
\begin{figure}
\center{
\includegraphics[scale=0.2]{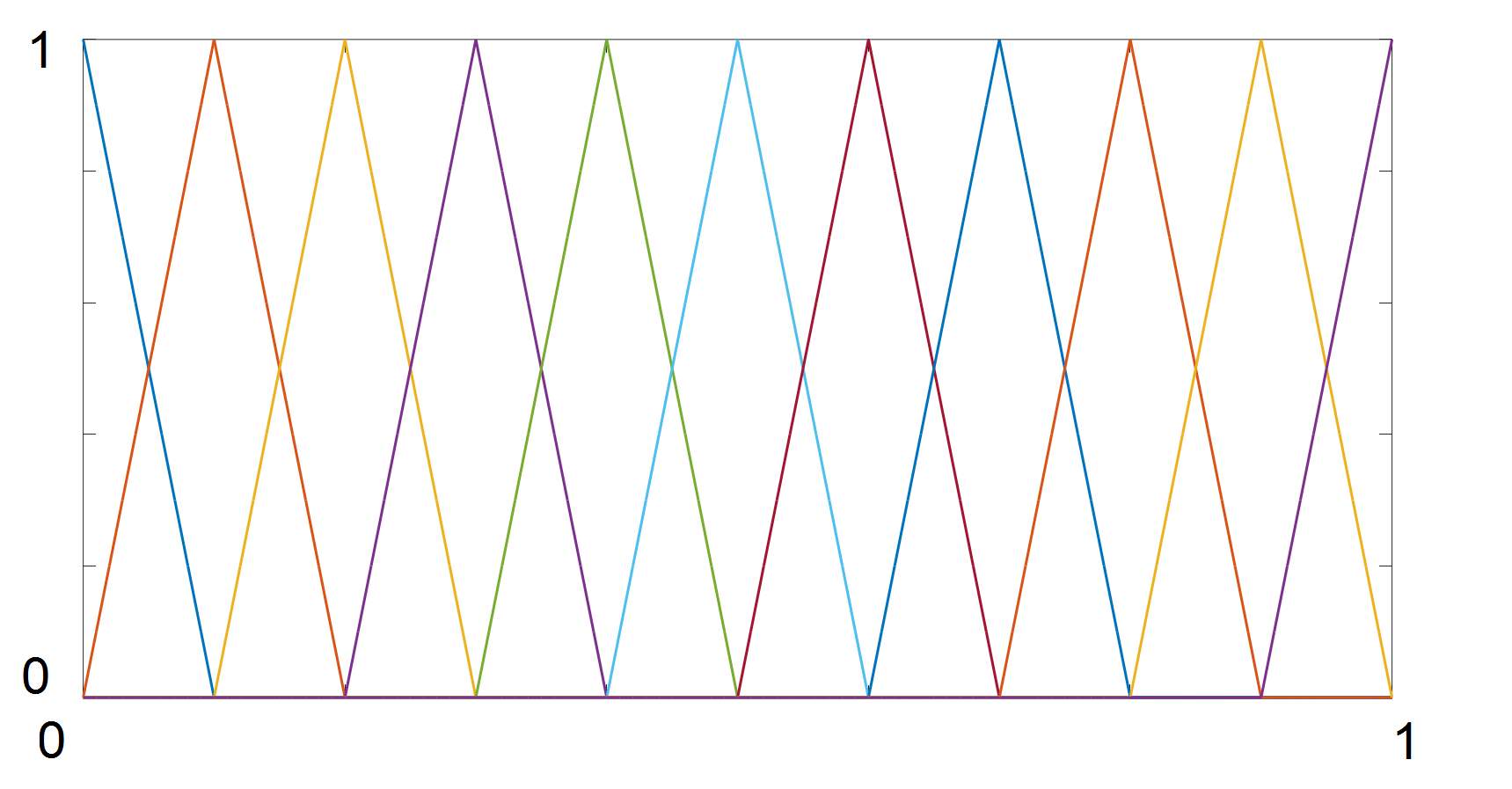}
\includegraphics[scale=0.2]{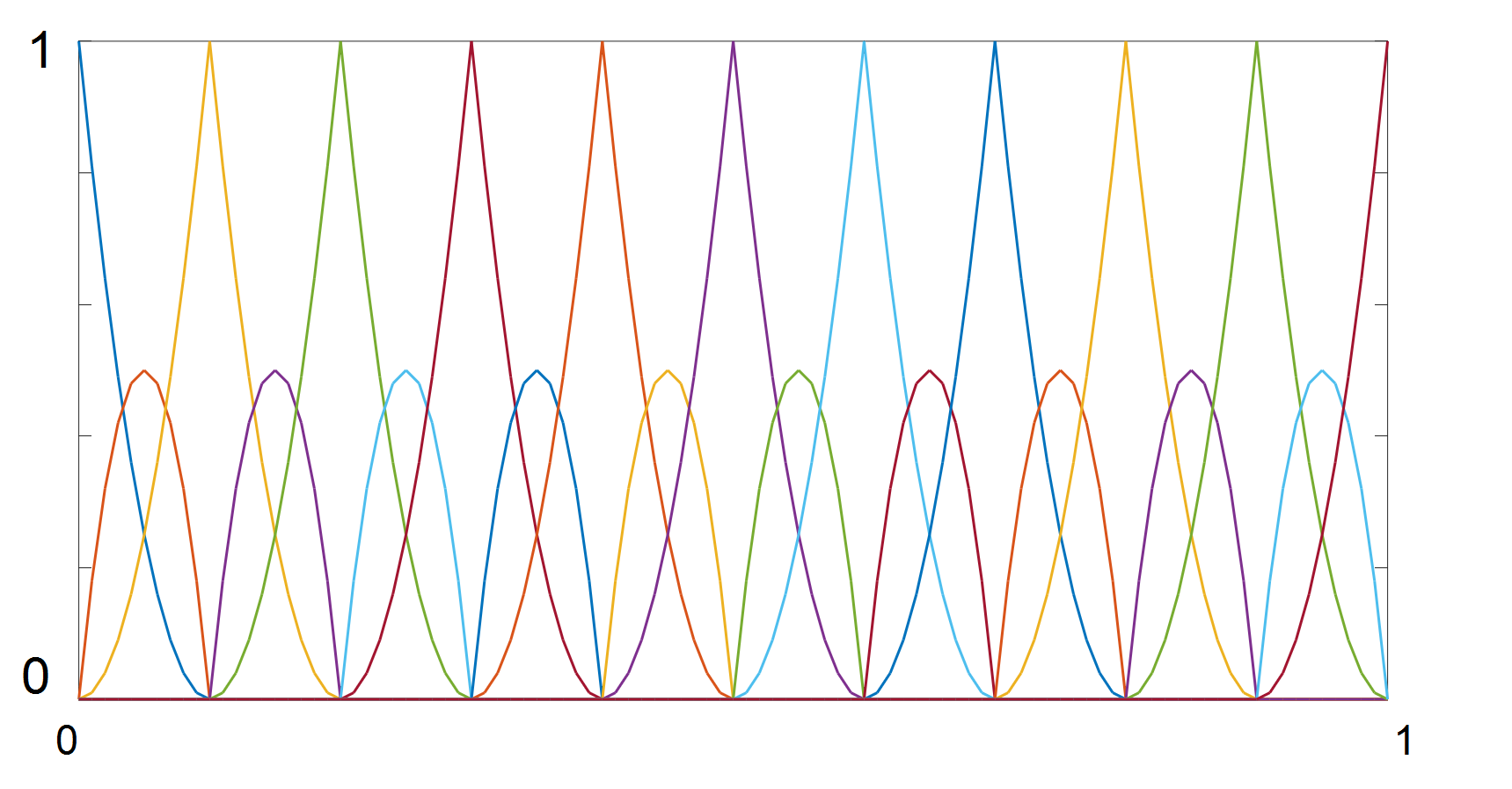}
}
\caption{Discrete spaces $U_h$ (left) and $V_h$ (right).}
\label{fig:Fig1}
\end{figure}

\begin{figure}
\center{
\includegraphics[scale=0.2]{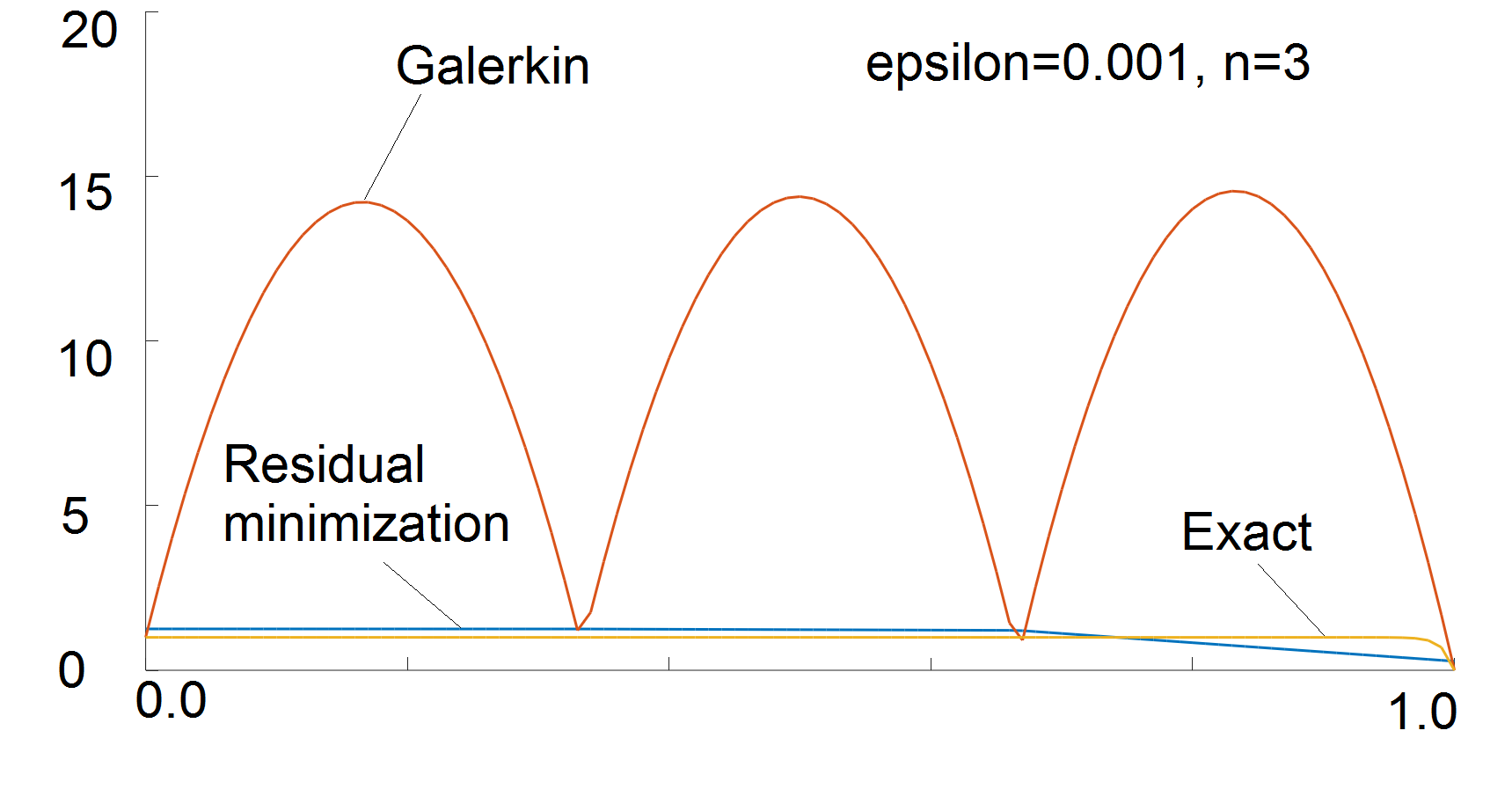}
\includegraphics[scale=0.2]{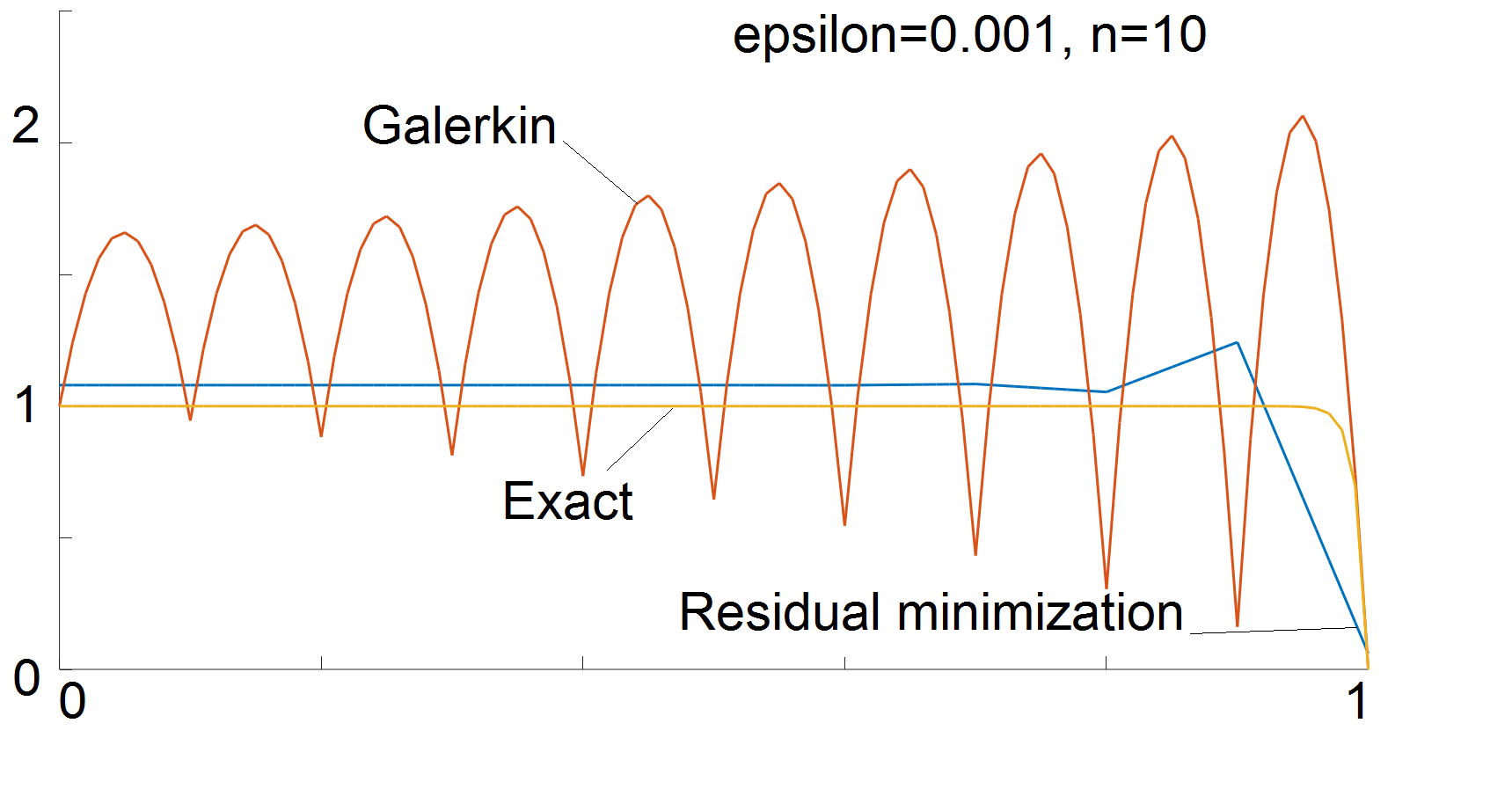} \\
\includegraphics[scale=0.2]{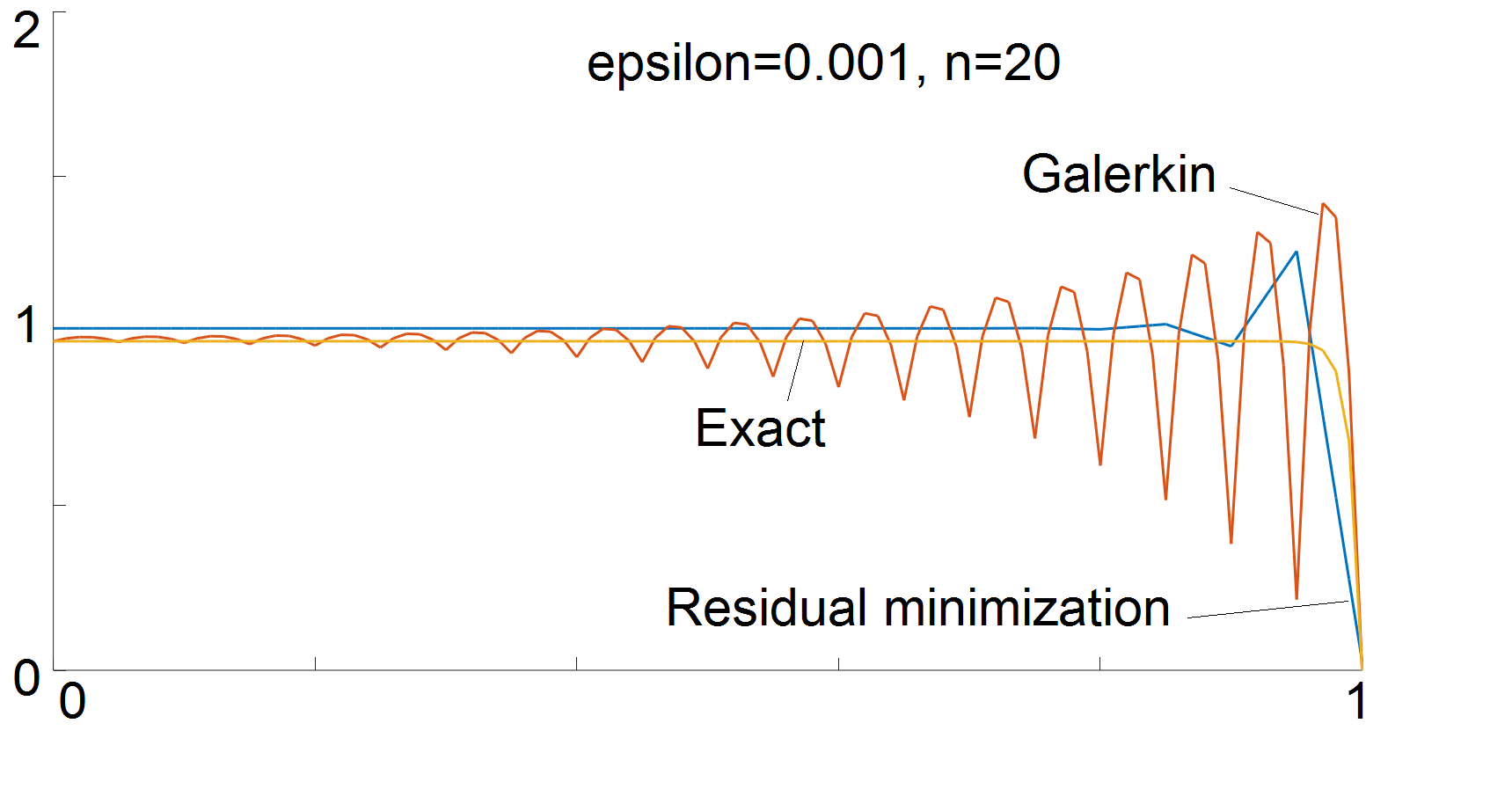}
\includegraphics[scale=0.2]{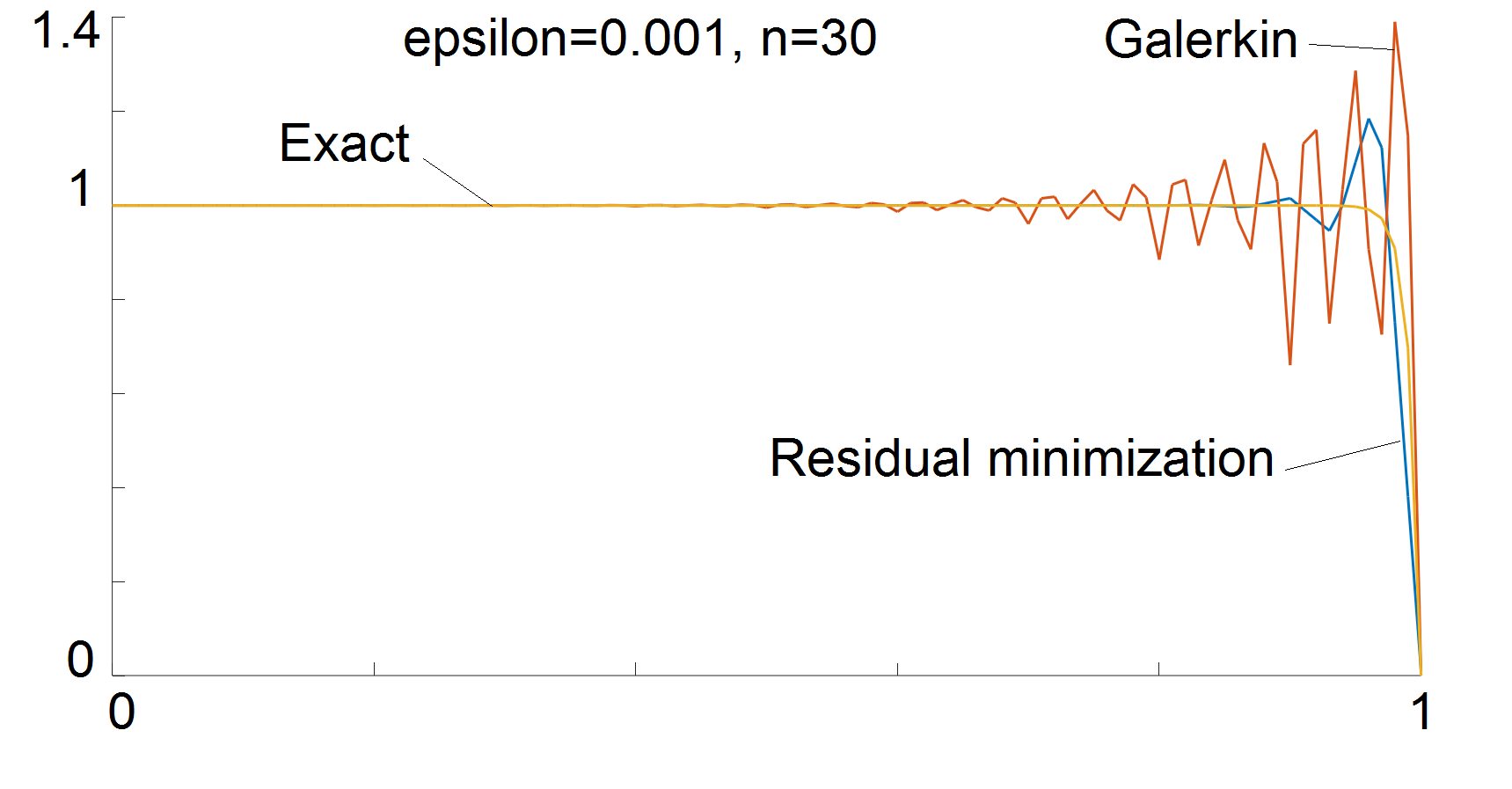} \\
\includegraphics[scale=0.2]{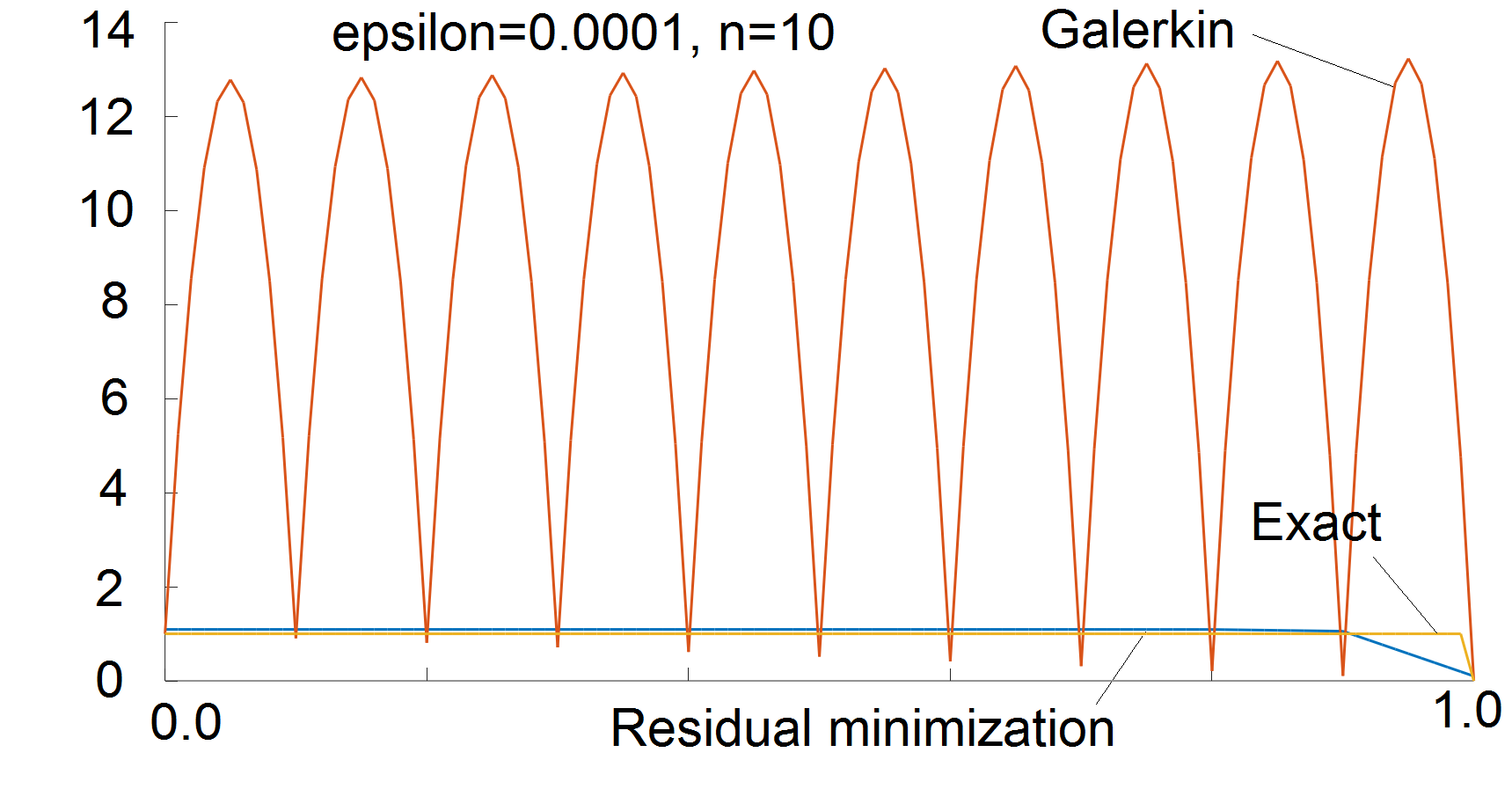} 
\includegraphics[scale=0.2]{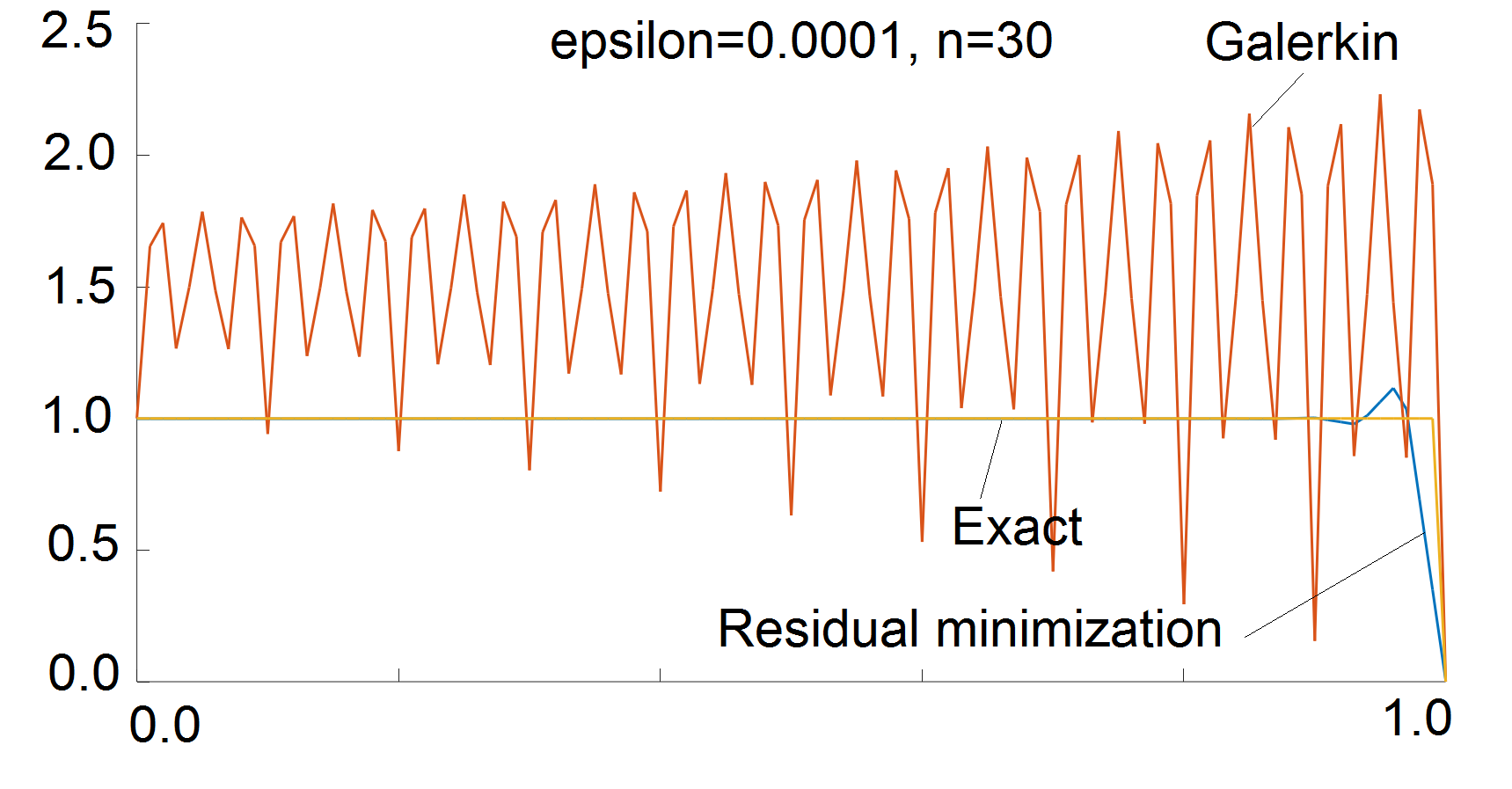} \\

}
\caption{Comparing the exact solution with the Galerkin method (where trial=test=quadratic B-splines with $C^0$ separators) and the residual minimization method (with linear B-splines for trial and quadratic B-splines with $C^0$ separators for test).}
\label{fig:Fig2}
\end{figure}

We define discrete spaces $U_h$ and $V_h$ as depicted in Figure~\ref{fig:Fig1}. Given a regular mesh, $U_h$ will be the space of piecewise linear and continuous functions; while $V_h$ will be the space of piecewise quadratics and continuous functions. From one side, we discretize 
formulation~\eqref{eq:varfor} using a standard Galerkin method where trial and test spaces are equal to $V_h$. On the other side, we discretize~\eqref{eq:varfor} by means of a residual minimization (RM) technique that uses $U_h$ as the trial space, and $V_h$ as the test space. Figure~\ref{fig:Fig2} compares the discrete solutions delivered by these two methods for different values of $\epsilon$, and different meshes parametrized by the number of elements $n$. We observe a superior performance of the residual minimization method, even though the approximation (trial) space used is poorer than that of the Galerkin method. 

\subsection{Outline of the paper}
\inblue{The structure of the paper is the following. Section 2 contains all the theoretical ingredients needed to understand our approach. Namely, Petrov-Galerkin formulations with optimal test functions (section 2.1); optimal test functions for an affine family of PDEs (section 2.2); the hierarchical compression of the optimal test functions matrix of coefficients (section 2.3); the fast hierarchical matrix-vector multiplication and related fast implementation of the GMRES solver (section 2.4); and the neural networks acceleration of the hierarchical matrix compression (section 2.5). Next, in section 3 we apply our automatic stabilization procedure to well-known unstable parametric model problems. First, for the two-dimensional Eriksson-Johnsson model problem (section 3.1); and next, for a two-dimensional Helmholtz equation (section 3.2). We write our conclusions in section 4. All the pseudo-code algorithms needed to follow our methodology have been shifted to the Appendix.}

\section{Theoretical ingredients}
\subsection{Petrov-Galerkin formulations with optimal test functions}
Let $U$ and $V$ be Hilbert spaces. We consider a general variational formulation of a PDE, which is to find $u\in U$ such that
\begin{equation}
b(u,v)=\ell(v)\,, \quad \forall v\in V\,, 
\label{eq:gen_weak}
\end{equation}
where $b: U\times V\to\mathbb R$ is a continuous bilinear form, and $\ell:V\to\mathbb R$ is a continuous linear functional (i.e., $\ell\in V'$, the dual space of $V$). 

The dual space $V'$ has a norm inherited by the norm of $V$.
Indeed, if $(\cdot,\cdot)_V$ denotes the inner-product of the Hilbert space $V$, then these norms are given by the following expressions:
$$
V\ni v\mapsto \|v\|_V:=\sqrt{(v,v)_V} \qquad\hbox{and}\qquad
V'\ni f\mapsto \|f\|_{V'}:=\sup_{\|v\|_V=1}f(v).
$$

We will assume well-posedness of problem~\eqref{eq:gen_weak}, which translates into the well-known inf-sup conditions (see, e.g.,~\cite[Theorem~2.6]{ern2013theory}):
\begin{subequations}
\begin{alignat}{2}
   & \exists \gamma>0 \hbox{ such that }\|b(w,\cdot)\|_{V'}\geq \gamma\|w\|_U\,, \quad\forall w\in U;\label{eq:bounded_below}\\
   & \{v\in V: b(w,v)=0\,, \forall w\in U\}=\{0\}.
    \end{alignat}
\end{subequations}
Notice that the continuity assumption on the bilinear form $b(\cdot,\cdot)$, also implies the existence of a constant $M\geq\gamma$ such that:
\begin{equation}\label{eq:continuity}
    \|b(w,\cdot)\|_{V'}\leq M \|w\|_U\,,\quad\forall w\in U.
\end{equation}


Given a discrete finite-element space $U_h\subset U$, a natural candidate to approximate the solution $u\in U$ to problem~\eqref{eq:gen_weak} is the residual minimizer 
\begin{equation}\label{eq:argmin}
 u_h = \argmin_{w_h \in U_h} \| b(w_h,\cdot) - \ell(\cdot) \|_{V'}\,.
\end{equation}
Indeed, combining~\eqref{eq:bounded_below},~\eqref{eq:argmin}, and~\eqref{eq:continuity}, the residual minimizer automatically satisfies the quasi-optimality property:
$$
\gamma\|u_h-u\|_U\leq 
\|b(u_h,\cdot)-\ell(\cdot)\|_{V'}=\inf_{w_h\in U_h}\|b(w_h,\cdot)-\ell(\cdot)\|_{V'}
\leq M \inf_{w_h\in U_h}\|w_h-u\|_U\,.
$$
Thus, the residual minimizer inherits stability properties from the continuous problem.

It is well-known (see, e.g.,~\cite{ChaEvaQiuCAMWA2014}) that the saddle-point formulation of residual minimization~\eqref{eq:argmin} becomes the mixed problem that aims to find $u_h\in U_h$, and a residual representative $r\in V$, such that:
\begin{subequations}\label{eq:mixed}
\begin{alignat}{3}
& (r,v)_V - b(u_h,v)  && = - \ell(v) , \quad  & \forall v\in V,\\
& b(w_h,r)  && = 0, \quad  & \forall w_h \in U_h.
\end{alignat}
\end{subequations}
Although it seems harmless, the mixed problem~\eqref{eq:mixed} is still infinite-dimensional in the test space $V$. To obtain a computable version, we introduce a discrete test space $V_h \in V$ that turns~\eqref{eq:mixed} into the fully discrete problem to find $u_h\in U_h$\footnote{Notice the abuse of notation. This new $u_h\in U_h$ solution of~\eqref{eq:mixed_discrete} does not equals the exact residual minimizer solution of~\eqref{eq:argmin}, or equivalently~\eqref{eq:mixed}.}, and a residual representative $r_h\in V_h$, such that:
\begin{subequations}\label{eq:mixed_discrete}
\begin{alignat}{3}
& (r_h,v_h)_{V} -  b(u_h,v_h)  &&= -\ell(v_h), \quad  & \forall v\in V_h,
\label{eq:mixed_discrete1}\\
& b(w_h,r_h) && = 0, \quad & \forall w_h \in U_h,
\label{eq:mixed_discrete2}
\end{alignat}
\end{subequations}
Problem~\eqref{eq:mixed_discrete} corresponds to the saddle-point formulation of a discrete-dual residual minimization, in which the dual norm $\|\cdot\|_{V'}$ in~\eqref{eq:argmin} is replaced by the discrete-dual norm $\|\cdot\|_{V_h'}$. Well-posedness and stability of~\eqref{eq:mixed_discrete} has been extensively studied in~\cite{MugVdZ_SINUM2020} and depends on a \emph{Fortin compatibility} condition between the discrete spaces $U_h$ and $V_h$. 
Moreover, once this condition is fulfilled (and in order to gain stability), it is possible to enrich
the test space $V_h$ without changing the trial space $U_h$. Obviously, this process will enlarge the linear system~\eqref{eq:mixed_discrete}. Nevertheless, we know that there is an equivalent linear system of the same size of the trial space, delivering the same $u_h\in U_h$ solving~\eqref{eq:mixed_discrete}. This is known as the Petrov-Galerkin method with optimal test functions, which we describe below. Our goal will be to make this generally impractical method practical.

Let us introduce now the concept of optimal test functions. For each $w_h\in U_h$, the optimal test function $Tw_h\in V_h$ is defined as the Riesz representative of the functional $b(w_h,\cdot)\in V_h'$, i.e., 
\begin{equation}\label{eq:opt_test}
(Tw_h,v_h)_V=b(w_h,v_h)\,, \quad \forall v_h\in V_h.
\end{equation}
Testing equation~\eqref{eq:mixed_discrete1} with optimal test functions $Tw_h$, using~\eqref{eq:opt_test} and~\eqref{eq:mixed_discrete2}, we arrive to the following Petrov-Galerkin system with optimal test functions:
\begin{equation}
\label{eq:PG}
\left\{
\begin{array}{ll}
\hbox{Find } u_h\in U_h \hbox{ such that}\\
b(u_h,Tw_h) = \ell(Tw_h)\,, & \forall w_h\in U_h\,.
\end{array}\right.
\end{equation}
In order to explicit a matrix expression for~\eqref{eq:PG}, let us set $U_h:=\operatorname{span}\{w_1,\dots,w_n\}$ and $V_h:=\operatorname{span}\{v_1,...,v_m\}$. 
Consider the matrix $\mathbb B$ linked to the bilinear form $b(\cdot,\cdot)$ such that its $(i,j)$-entry is $\mathbb B_{ij}=b(w_j,v_i)$.
Analogously, we consider the Gram matrix $\mathbb G$ linked to the inner product $(\cdot,\cdot)_V$ such that $\mathbb G_{ki}=(v_k,v_i)_V$. 
The optimal test space is defined as $V_h^{\hbox{\tiny{opt}}}:=
\operatorname{span}\{Tw_1,...,Tw_n\}$.
Thus, using~\eqref{eq:opt_test}, we observe that the matrix containing the coefficients of optimal test functions when expanded in the basis of $V_h$ is $\mathbb W:=\mathbb G^{-1}\mathbb B$. Moreover, the Petrov-Galerkin system~\eqref{eq:PG} becomes
\begin{equation}
\label{eq:PG_sys}
\mathbb B^T\mathbb W\,x = 
(L^T\mathbb W)^T,
\end{equation}
where the vector $x$ contains the coefficients of the expansion of $u_h$ in the basis of $U_h$, $L^T=[\ell(v_1)\,\cdots\, \ell(v_m)]$, and we have used the fact that $\mathbb G^T=\mathbb G$. Therefore, if we aim to solve the Petrov-Galerkin linear system~\eqref{eq:PG_sys}, an optimized matrix-vector multiplication to perform $\mathbb W\,x$ and $L^T\mathbb W$ becomes critical. Section~\ref{sec:compression} is devoted to study the hierarchical compression of $\mathbb W$, which allows for fast vector-matrix multiplications.

\subsection{Optimal test functions for an affine family of parametric PDEs}
\label{sec:affine}
Assume we want to solve parametric PDEs in variational form, i.e., 
$$
\hbox{Find }\quad u_\mu\in U\quad \hbox{ such that }\quad
b_\mu(u_\mu,v)=\ell_\mu(v), \quad\forall v\in V,
$$
where for each set of parameters $\mu\in \mathcal P\subset \mathbb R^d$, the bilinear form $b_\mu(\cdot,\cdot)$ is continuous and inf-sup stable, with constants that may depend on $\mu$. Moreover, we assume that  $b_\mu(\cdot,\cdot)$ has the affine decomposition:
$$
b_\mu(\cdot,\cdot) =b_0(\cdot,\cdot)+\sum_{k=1}^d\theta_k(\mu)b_k(\cdot,\cdot)\,,
$$
where $\theta_k : \mathcal P \to \mathbb R$ and $b_k:U\times V\to\mathbb R$ are accessible and easy to compute.
When the trial and test spaces are discretized, the bilinear form $b_\mu(\cdot,\cdot)$ induces a matrix of the form:
$$
\mathbb B_\mu= \mathbb B_0 +\sum_{k=1}^d\theta_k(\mu)\mathbb B_k\,.
$$
Thus, the matrix of coefficients of optimal test functions $\mathbb W_\mu:=\mathbb G^{-1}\mathbb B_\mu$ becomes in this case:
\begin{equation}\label{eq:W}
\mathbb W_\mu = \mathbb G^{-1}\mathbb B_0 +\sum_{k=1}^d\theta_k(\mu)\mathbb G^{-1}\mathbb B_k\,.
\end{equation}
\inblue{Equation~\eqref{eq:epsilonPG} shows the particular form 
that expression~\eqref{eq:W} gets for the Eriksson-Johnsson model problem, where knowing $\mathbb G^{-1}\mathbb B_0$ and $\mathbb G^{-1}\mathbb B_1$ implies the knowledge of $\mathbb W_\epsilon$ for any $\epsilon>0$.
}

\subsection{Hierarchical compression of the optimal test functions matrix of coefficients}\label{sec:compression}

The hierarchical matrices has been introduced by Hackbush \cite{Hackbush}.
The main idea of the hierachical compression of a matrix is to store the matrix in a tree-like structure, where:
\begin{itemize}
    \item the root node corresponds to the whole matrix;
    \item the root node has some number of sons (in our approach 4 sons) corresponding to submatrices of the main matrix;
    \item each node can have sons (in our approach 4 sons) corresponding to submatrices (blocks), or can be a leaf representing the corresponding matrix (block);
    \item each leaf stores its associated matrix in the SVD compressed form \ingreen{or as zero matrix};
    \item at each node, the decision about storing the block in SVD form, or either dividing the block into submatrices, depends on \ingreen{an} admissibility condition of the block.
\end{itemize}
\begin{figure}
\begin{center}
\includegraphics[width=0.6\textwidth] {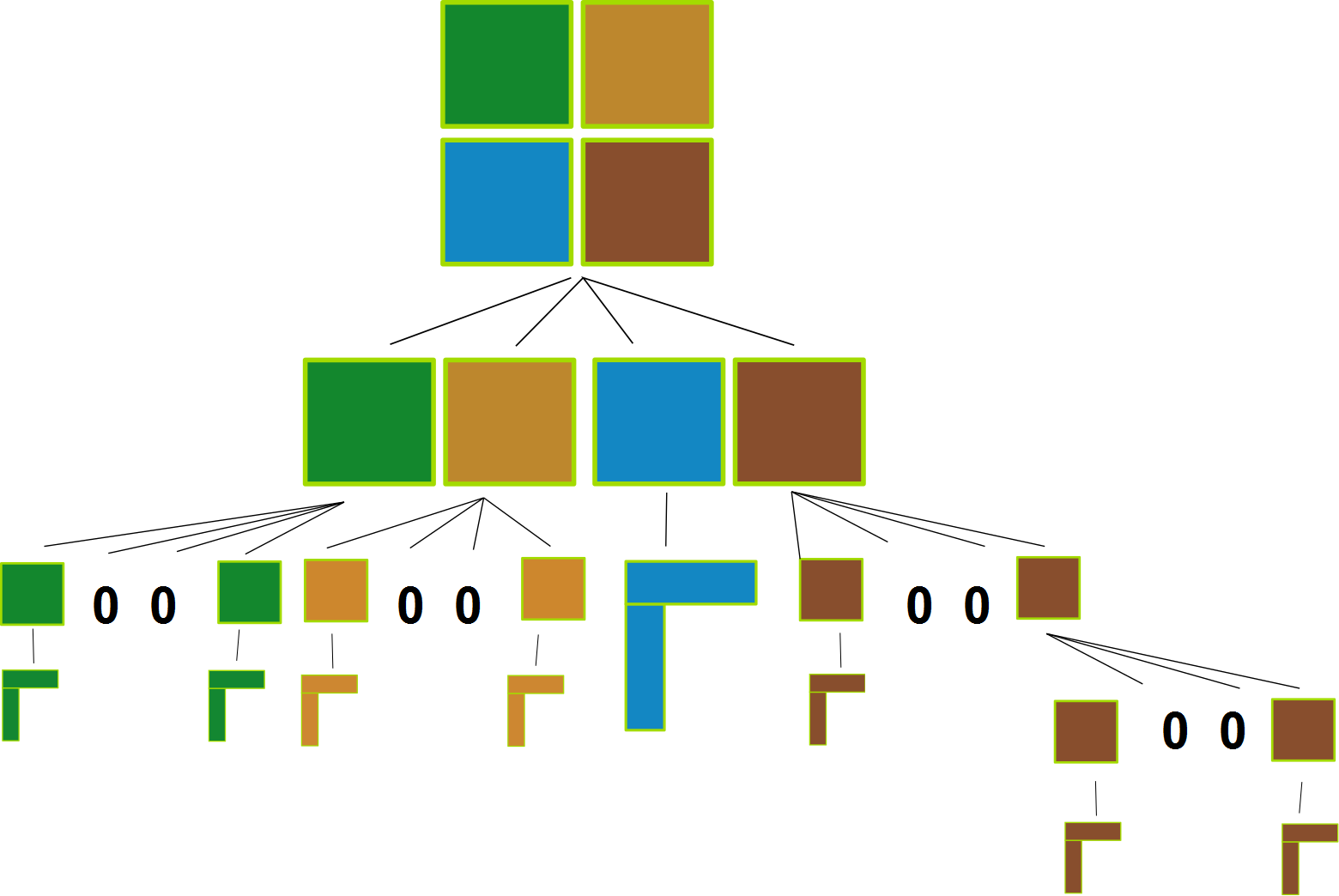}
 \caption{Hierarchical compression of a matrix}
 \label{fig:tree}
\end{center} 
 \end{figure}
Exemplary hierarchical compression of the matrix in a form of a tree is presented in Figure~\ref{fig:tree}; 
%
while the algorithm for compression of the matrix into the hierarchical matrix format is presented in Algorithm \ref{Alg1}.

The admissibility condition controls the process of creation of the tree, it  allows to decide if the matrix should be divided (or not) into submatrices.
In our case, the admissibility condition is defined by 
\begin{enumerate}
    \item The size of the matrix: if the matrix is bigger than a pre-defined maximal admissible size $l> >1$, then the matrix should be divided into submatrices;
    \item The first $r$ singular values: if the $r+1$ singular value is greater than a pre-defined threshold $\delta>0$, then the matrix should be divided into submatrices.
\end{enumerate}

\begin{figure}
\begin{center}
  \includegraphics[width=0.4\textwidth]{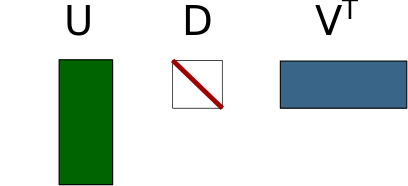}
\caption{Reduced Singular Value Decomposition.}
\label{fig:rSVD}
\end{center}
\end{figure}

In the leaves of the tree, we perform a reduced Signular Value Decomposition (rSVD). A \ingreen{reduced singular value decomposition of} a $(n\times m)$-matrix $\mathbb M$ of rank $k$ is a factorisation of the form 
$$
\mathbb M=\mathbb U\,\mathbb D\,\mathbb V^T,
$$
with unitary matrices $\mathbb U \in \mathbb{R}^{n\times \ingreen{k}}$ and $\mathbb V \in \mathbb{R}^{m\times \ingreen{k}}$, and a diagonal matrix $\mathbb D \in \mathbb{R}^{\ingreen{k}\times \ingreen{k}}$ where the diagonal entries are $\mathbb D_{11} \ge \mathbb D_{22} \ge\cdots\ge \mathbb D_{kk}>0$. (see Figure \ref{fig:rSVD}).
The diagonal entries of $\mathbb D$ are called the singular values of $\mathbb M$.
The computational complexity of the reduced SVD is ${\cal O}((m+n)\ingreen{k}^2)$.


\subsection{Matrix-vector multiplication with $\mathbb H$-matrices and GMRES solver speedup}
\begin{figure}
    \centering
 \includegraphics[width=0.5\textwidth]{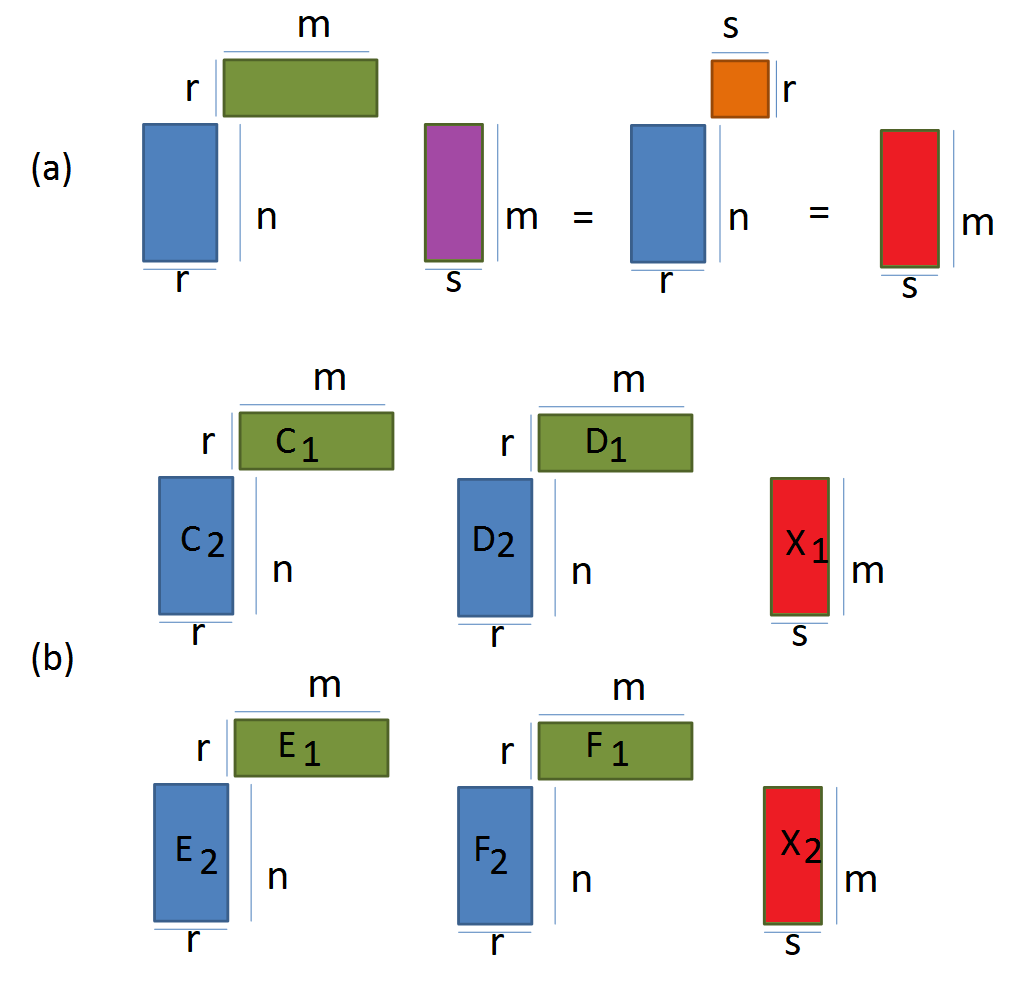}
    \caption{(a) Computational cost of matrix-vector multiplication with compressed $\mathbb H$-matrix and $s$ vectors ${\cal O}(rms+rns)$ when $n=N>>r$ reduces to ${\cal O}(Nrs)$. (b) Multiplication of a matrix compressed into four SVD blocks by the vector partitioned into two blocks, following  $\begin{bmatrix} \ingreen{C_2*(C_1}*X_1)+\ingreen{D_2*(D_1}*X_2) \\  \ingreen{E_2*(E_1}*X_1)+\ingreen{F_2*(F_1}*X_2) \end{bmatrix}$. The resulting computational cost is ${\cal O}(Nrs)$.}
\label{fig:MVMult}
\end{figure}

The computational cost of matrix-vector multiplication using a compressed $\mathbb H$-matrix \ingreen{of rank $r$} and $s$ right-hand side vectors is ${\cal O}((m+n)rs)$. \inblue{This is illustrated in Figure ~\ref{fig:MVMult}(a).}

The multiplication of a matrix compressed into SVD blocks is performed recursively as illustrated in \inblue{Figure~\ref{fig:MVMult}(b).}
The resulting computational cost of the multiplication is ${\cal O}(Nrs)$, \ingreen{where $N:=\max\{n,m\}$}.

The GMRES algorithm employed for computing the solution, includes multiplications of the problem matrix by vectors (see line 1, line 4, and line 5 in Algorithm \ref{Alg4}). The application of the stabilization matrix requires replacement of the $\mathbb B$ matrix by $\mathbb B^T\mathbb H$. 
If we can apply matrix-vector multiplication $\mathbb H\, x$ in a linear cost, we can say that our stabilization comes \emph{for free}.


\subsection{Neural network learning the hierarchical matrices}
\begin{figure}
    \centering
 \includegraphics[width=0.9\textwidth]{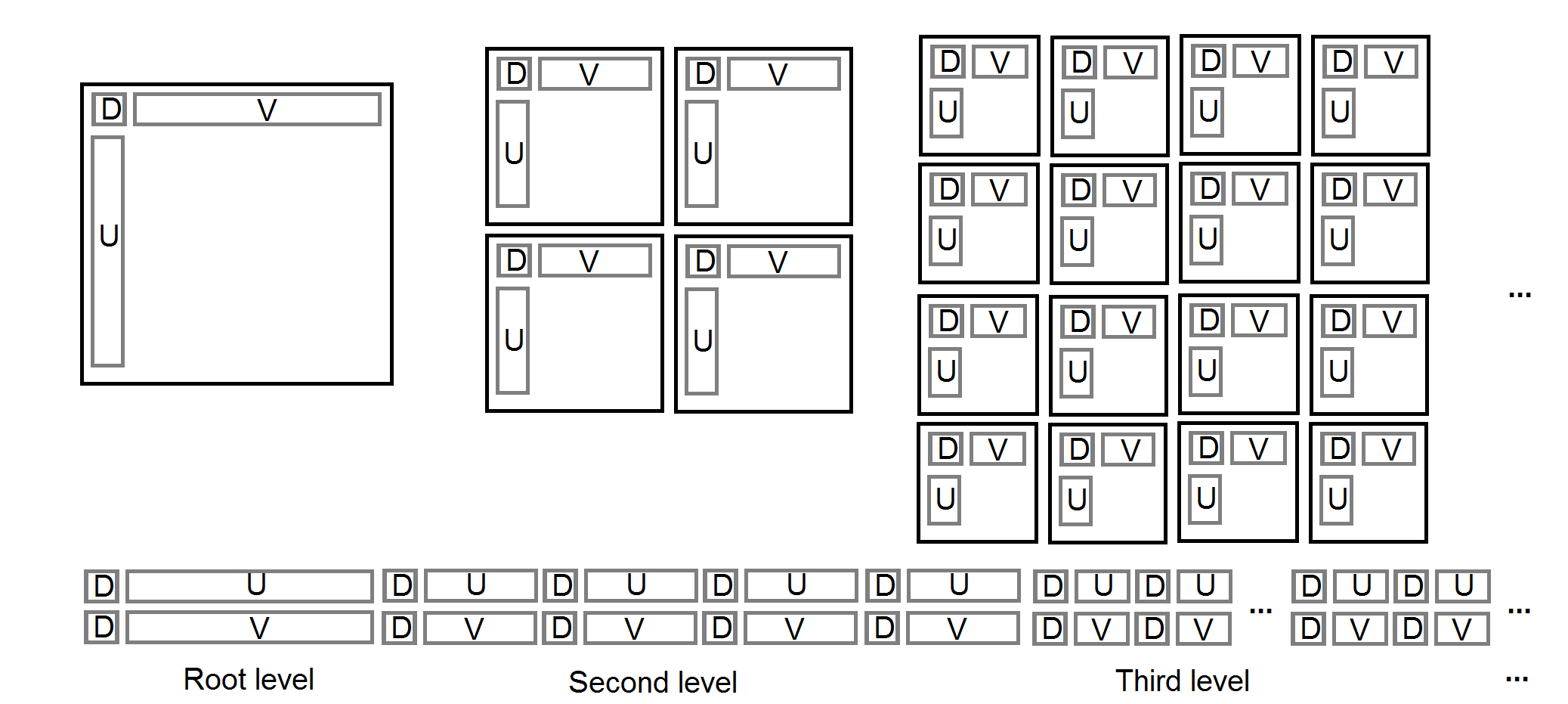}
  \caption{Compression of truncated SVD decompositions of matrices, starting from the root level, second level, third level, and the following levels. }
\label{fig:compression}
\end{figure}

The hierarchical matrix is obtained by constructing a tree with SVD decompositions of different blocks of the full matrix. The root level corresponds to the entire matrix, and  the children correspond to sub-blocks. 
Only the leaf nodes have blocks stored in the SVD decomposition format. The most expensive part of the compression algorithm is checking the admissibility condition. In particular, checking if a given block has $r$ singular values smaller than $\delta$, and whether we partition or run SVD decomposition.
The SVD
data for the blocks of different sizes can be precomputed and stored in a list, see Figure~\ref{fig:compression}.
From the set $\mathcal P\subset \mathbb R^d$ of PDE-parameters, we can construct the neural network
\begin{equation}
\mathcal P \ni \mu \rightarrow \operatorname{DNN}(\mu)
=\{\mathbb U_i(\mu),\mathbb D_i(\mu),\mathbb V_i(\mu) \}_{i=1,...,N_B}
\end{equation}
where $\operatorname{DNN}(\mu)$ is the list of SVD decompositions for all $N_B$ blocks of different dimensions. 

\begin{figure}
    \centering
 \includegraphics[width=0.9\textwidth]{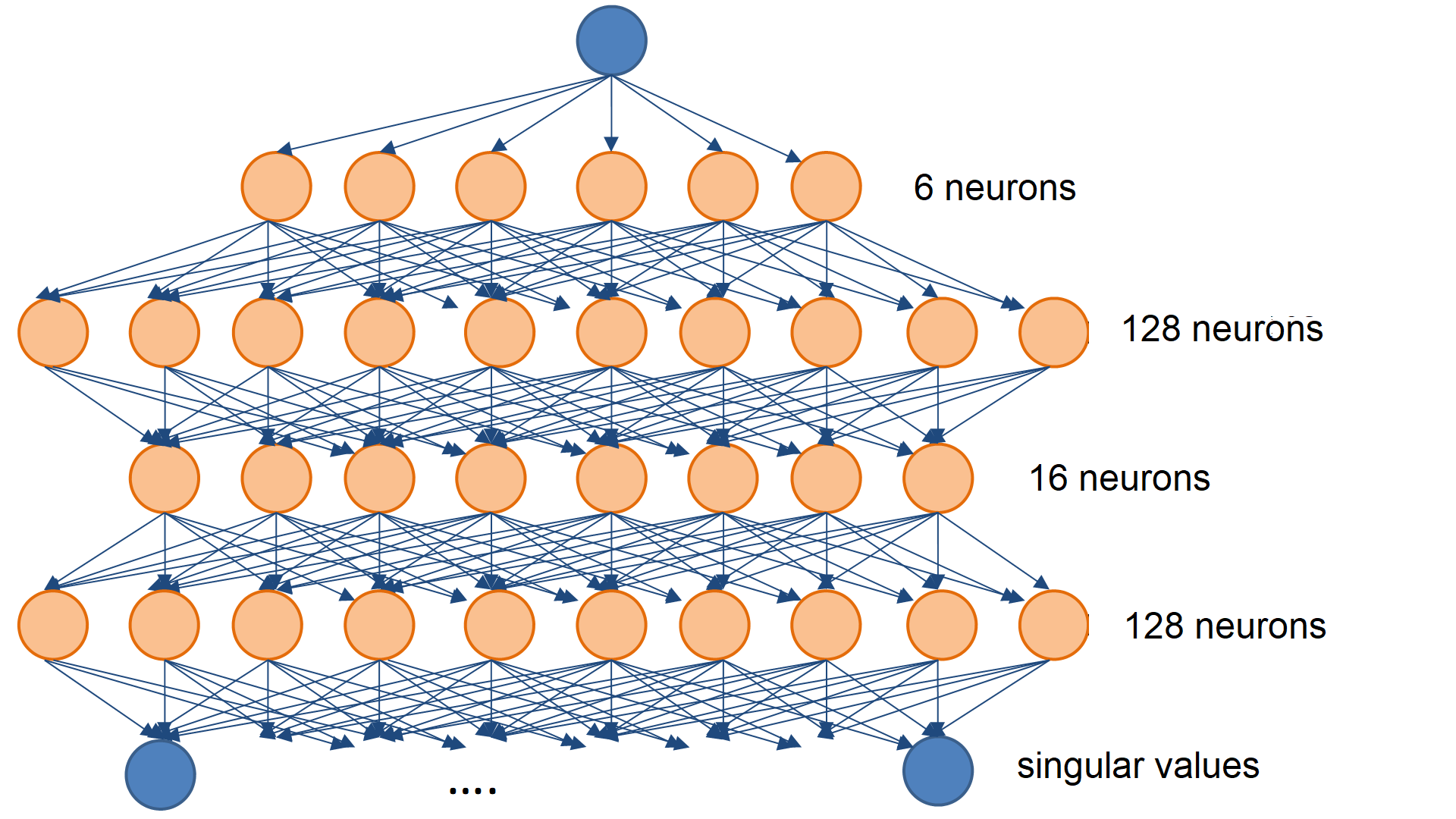}
    \caption{The architecture of the neural network learning the singular values for blocks of the matrix. }
\label{fig:training16x16}
\end{figure}

Unfortunately, training the neural networks $\mathbb U_i(\mu)$ and $\mathbb V_i(\mu)$ for different blocks, as functions of $\mu$ does not work. However, it is possible to train the singular values for different blocks as function of $\mu$, i.e., 
\begin{equation}
    \mathcal P\ni \mu \rightarrow \operatorname{DNN}(\mu)
    =\{\mathbb D_i(\mu)\}_{i=1,...,N_B}\,.
\end{equation}
Knowing $\mathbb D_i(\mu)$ a priori for a given $\mu$ allows us to construct the structure of the hierarchical matrix, and call truncated SVD only for leaves of the matrix, which considerably speeds up the compression algorithm. An illustration of the architecture of the neural network used to learn singular values is depicted in Figure~\ref{fig:training16x16}.
%


\section{Application}
\subsection{Two-dimensional Eriksson-Johnson problem}
Given $\Omega=(0,1)^2\subset\mathbb R^2$ and $\beta=(1,0)$, we seek the solution of the advection-diffusion problem
\begin{equation}
\left\{
\begin{array}{rl}
-\epsilon\,\Delta u + \beta\cdot\nabla u=0 & \hbox{in }\Omega\\
u=\sin(k\pi y)\chi_{\{x=0\}} & \hbox{over }\partial\Omega\,,
\end{array}\right.
\label{eq:Erikkson}
\end{equation}
where $\chi_{\{x=0\}}$ denotes the characteristic function over the the inflow boundary $x=0$.
\begin{figure}[h]
\centering
\includegraphics[scale=0.28]{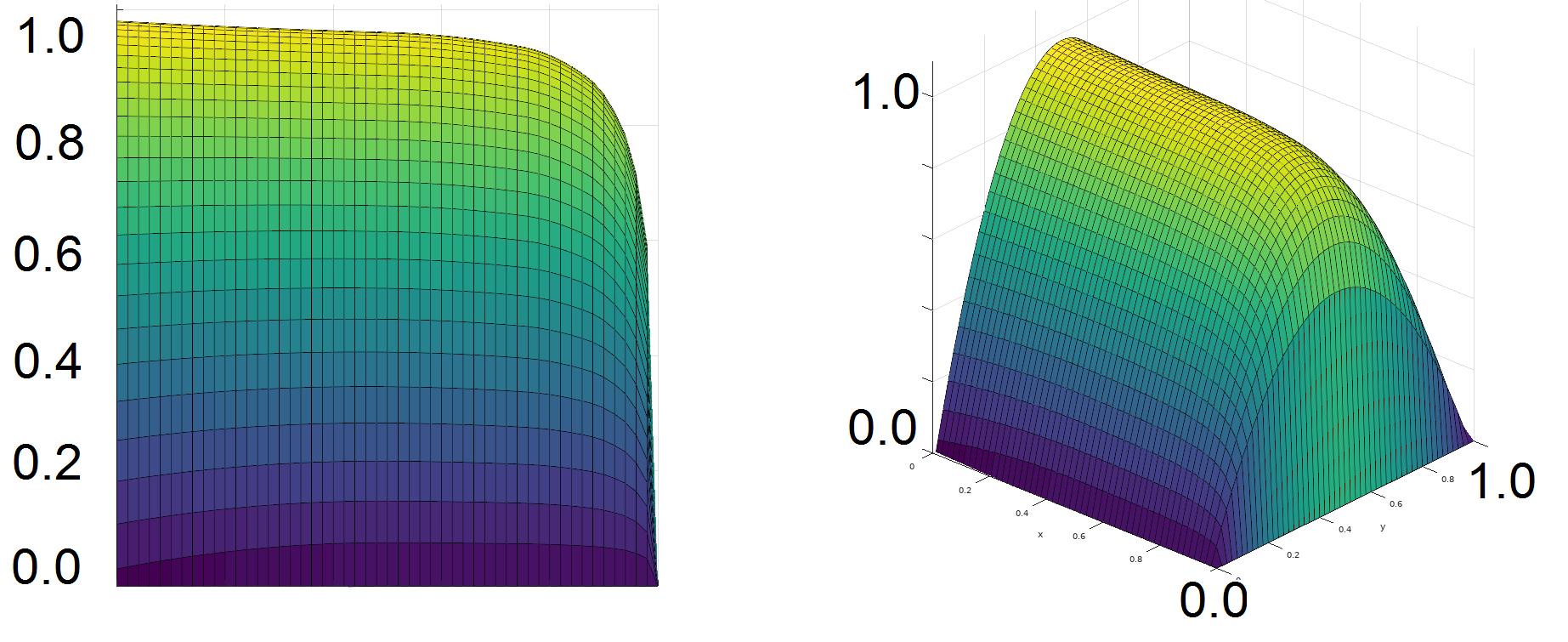}
\caption{Boundary layer of Eriksson-Johnson problem for $\epsilon=0.01$.}
\label{fig:Erikkson_problem}
\end{figure}
The problem is driven by the inflow Dirichlet boundary condition and develops a boundary layer of width $\epsilon$ near the outflow $x = 1$, as shown in Figure~\ref{fig:Erikkson_problem}.


The weak form with a general Dirichlet boundary data $g\in H^{1\over2}(\partial\Omega)$ will be to find 
$u^\epsilon\in H^1(\Omega)$ such that 
\begin{equation}\label{eq:weakform}
\left\{
\begin{array}{rl}
\underbrace{\displaystyle\epsilon\int_\Omega \nabla u^\epsilon\cdot \nabla v 
+\int_\Omega(\beta\cdot \nabla u^\epsilon) v}_{\epsilon \,b_1(u^\epsilon,v)+b_0(u^\epsilon,v)=:b_\epsilon(u^\epsilon,v)} = 0 , & \forall v\in H_0^1(\Omega),\\
u^\epsilon = g , & \hbox{over }\partial\Omega.
\end{array}
\right.
\end{equation}
To simplify the discussion, we approximate the solution as tensor products of one-dimensional B-splines basis functions $\{B_{i;p}(x)B_{j;p}(y)\}_{i,j}$ of uniform order $p$ in all directions. This discrete trial space $U_h^p\subset H^1(\Omega)$ will be split as $U_h^p=U_{h,0}^p + U^p_{h,\partial\Omega}$, where $U^p_{h,0}\subset H_0^1$ contains all the basis functions vanishing at $\partial\Omega$; and $U^p_{h,\partial\Omega}$ is the complementary subspace containing the basis functions associated with boundary nodes. Our discrete solution will be $u_h^\epsilon=u_{h,0}^\epsilon+u^\epsilon_{h,g}$, where $u^\epsilon_{h,0}\in U_{h,0}^p$ is unknown and $u^\epsilon_{h,g}\in U_{h,\partial\Omega}^p$ is directly obtained using the boundary data $g$. 

We build the test space using a larger polynomial order $q>p$. That is, we use the tensor product of one-dimensional B-splines basis functions $\{B_{s;q}(x)B_{t;q}(y)\}_{s,t}$ of order $q$ and vanishing over $\partial\Omega$. This discrete test space will be denoted by $V^q_{h,0}\subset H_0^1(\Omega)$.
The discrete residual minimization problem will be to find $u^\epsilon_{h,0}\in U_{h,0}^p$ and $r_h\in V_{h,0}^q$ such that
\begin{subequations}\label{eq:epsilon_mixed}
\begin{alignat}{3}
&  (\nabla r_h,\nabla v_h)_{L^2(\Omega)} - b_\epsilon(u^\epsilon_{h,0},v_h) && = b_\epsilon(u^\epsilon_{h,g},v_h)\,, && \quad\forall v_h\in V_{h,0}^q\,;\\
& b_\epsilon(w_{h},r_h) && = 0 \,, && \quad\forall w_{h}\in U_{h,0}^p\,.
\end{alignat}
\end{subequations}

The reduced matrix system associated with~\eqref{eq:epsilon_mixed} takes the form:
\begin{equation}\label{eq:epsilonPG}
\mathbb B_\epsilon^T\mathbb W_\epsilon x= (L^T\mathbb W_\epsilon)^T, 
\quad\hbox{ where } 
\mathbb W_\epsilon =\mathbb G^{-1} \mathbb B_\epsilon = \epsilon\,\mathbb G^{-1}\mathbb B_1 + \mathbb G^{-1}\mathbb B_0. 
\end{equation}

\begin{figure}
    \centering
 \includegraphics[width=0.5\textwidth]{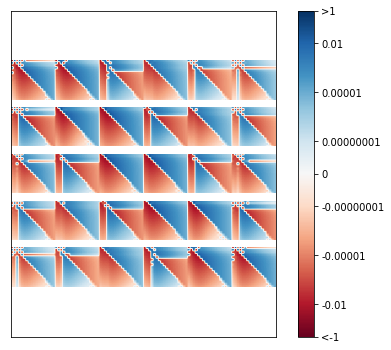}
    
    \caption{The entire matrix $\mathbb W_\epsilon$ of the coefficients of the optimal test functions as the input of the first hierarchical level. }
\label{fig:partition1x1}
\end{figure}
\begin{figure}
    \centering
\includegraphics[width=0.495\textwidth]{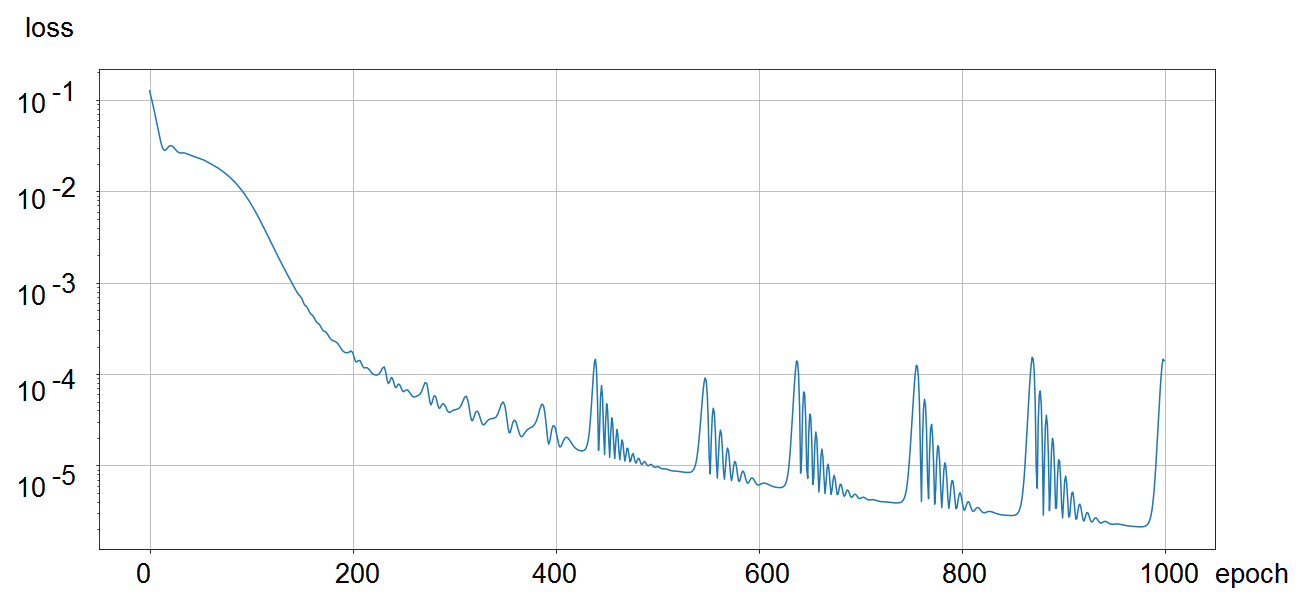}
\includegraphics[width=0.495\textwidth]{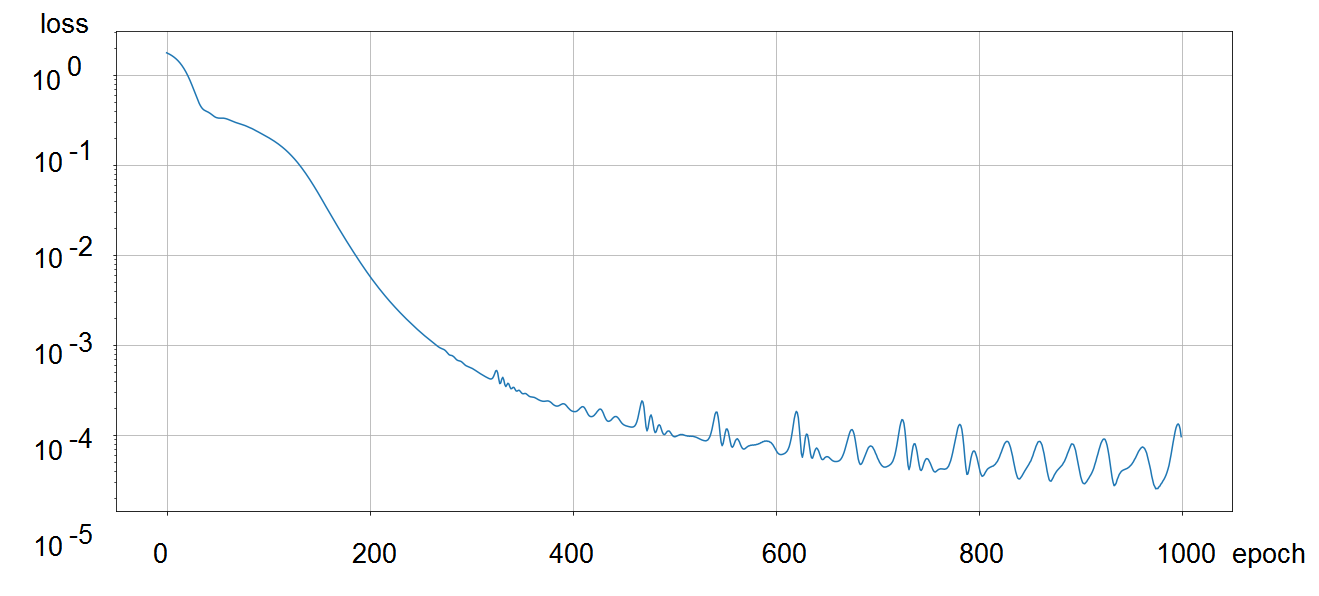}
\includegraphics[width=0.495\textwidth]{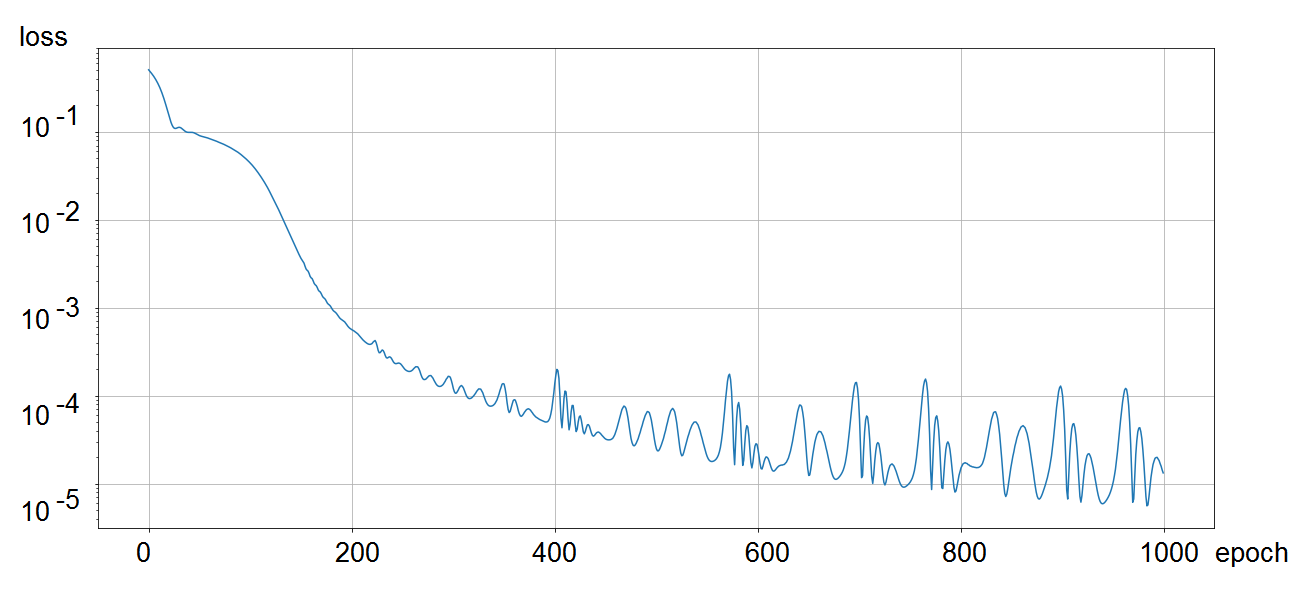}
\includegraphics[width=0.495\textwidth]{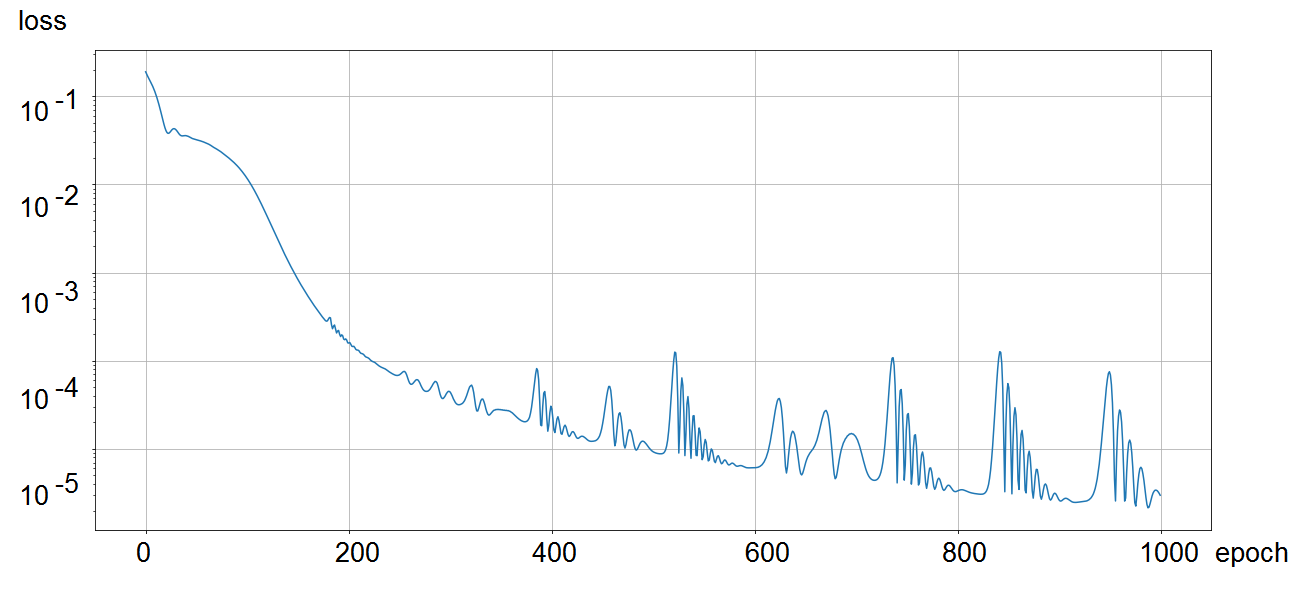}
    \caption{Erikkson-Johnsson problem. Training of the neural network for the diagonal matrices $\mathbb D_{ij}(\epsilon)$ at several levels.}
\label{fig:training1x1}
\end{figure}

We train a neural network for the diagonal $\mathbb D(\epsilon)$ matrix of the SVD decomposition $[\mathbb U(\epsilon),\mathbb D(\epsilon),\mathbb V(\epsilon)]$ of the entire matrix $\mathbb W_\epsilon$ (see Figure~\ref{fig:partition1x1}), as well as
$\mathbb D_{ij}(\epsilon)$ from the SVD decompositions $\{\mathbb U_{ij}(\epsilon),\mathbb D_{ij}(\epsilon),\mathbb V_{ij}(\epsilon)\}_{i=1,...,j^2}$ of sub-matrices obtained by $j \times j$ partitions of $\mathbb W_\epsilon$, for $j=2,4,8,16$. 
The white parts correspond to the boundary nodes, where we have enforced the boundary conditions.
The convergence of the training procedures is presented in Figures \ref{fig:training1x1}.
%
%











\begin{figure}
    \centering
 \includegraphics[width=0.45\textwidth]{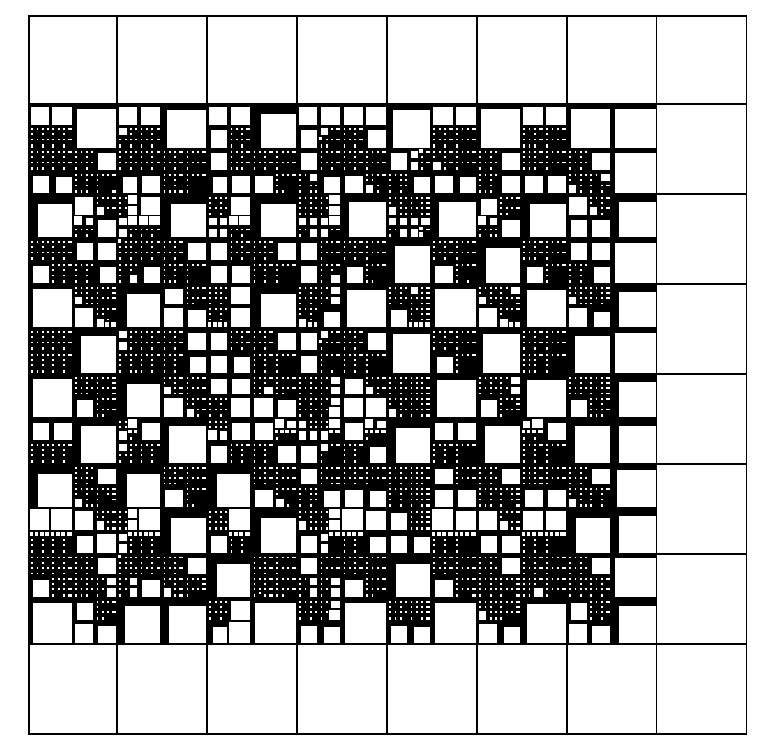}
 \includegraphics[width=0.45\textwidth]{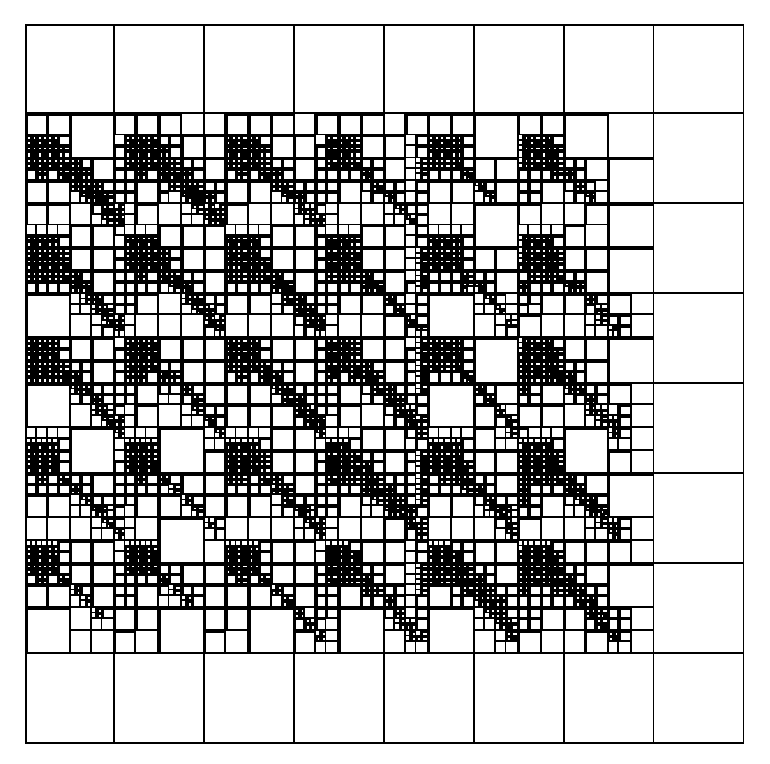}
    \caption{Matrix of the coefficients of the optimal test functions compressed by using recursive SVD algorithm for $\epsilon=0.1$ and $\epsilon=0.000001$. The compression algorithm removes singular values smaller than $\delta=10^{-7}$, and employs $r=16$.}
\label{fig:compression1}
\end{figure}

In an online stage, for a given diffusion coefficient $\epsilon$, we perform the compression of the matrix $\mathbb W_\epsilon$ 
into the hierarchical matrix $\mathbb H_\epsilon$ using Algorithm 6, where the admissibility condition is now provided by the trained neural network.
%
The compressed hierarchical matrices 
are illustrated in Figure~\ref{fig:compression1}.

Having the compressed matrix $\mathbb H_\epsilon$,
%
we employ the GMRES algorithm~\cite{Saad} for solving the linear system~\eqref{eq:epsilonPG}.
To avoid the computation with a dense $\mathbb B_\epsilon^T\mathbb H_\epsilon$ matrix,
we note that the GMRES method involves computations of the residual 
$R = \mathbb B_\epsilon^T\mathbb H_\epsilon\, x -(L^T\mathbb H_\epsilon)^T$
and the hierarchical matrix $\mathbb H_\epsilon$ enables matrix-vector multiplications of $\mathbb H_\epsilon\, x$ and $L^T\mathbb H_\epsilon$ in a quasi-linear computational cost.
%
%

\begin{table}
\begin{tabular}{ccccccc}
$\epsilon$ & {\bf Compress } & {\bf Compress } & {\bf ${\mathbb H}_{\epsilon}*x$} & {\bf ${\mathbb A}*{\mathbb H}_{\epsilon}x$} & {\bf \# iter} & {\bf Total} \\
& {\bf ${\mathbb H}_{\epsilon}$ flops} & {\bf ${\mathbb H}_{\epsilon}$ flops} & {\bf flops} & {\bf flops} & \inblue{\bf GMRES } & {\bf flops}  \\
& {\bf with DNN} & {\bf without DNN} &  & & \inblue{{\bf ${\mathbb H}$-matrix}} & \\
 \hline
0.1 & 14,334 & 158,946 & 41,163 & 9,152 & 90 & 3,728,166 \\ 
$10^{-6}$ & 31,880 & 117,922 & 33,100 & 9,002 & 76 & 3,232,392  \\
\hline
\vspace{5pt}
\end{tabular}
\caption{Computational costs of the stabilized Erikkson-Johnsson solver  using neural networks, hierarchical matrices and GMRES solver.}
\label{tab:input_2D}
\end{table}

\begin{table}
\begin{tabular}{cccc}
$\epsilon$ &  {\bf \# iter GMRES Galerkin} & {\bf Flops per iteration} & {\bf Total flops}  \\
 \hline
0.1 &  89 & 17,536 & 1,560,704\\ 
$10^{-6}$ & 65 & 17,536 & 1,139,840 \\
\hline
\vspace{5pt}
\end{tabular}
\caption{Computational costs of Galerkin method for the Erikkson-Johnsson problem using GMRES solver.}
\label{tab:input_2Da}
\end{table}

\inred{In Table \ref{tab:input_2D} we summarize the computational costs of our solver for two values of $\epsilon=\{0.1,0.000001\}$. The computational mesh was a tensor product of quadratic B-splines with 26 elements along $x$-axis and quadratic B-splines with 10 elements along $y$-axis. As we can read from the second and third columns, the DNN speeds up the compression process of the matrix of optimal test function's coefficients around ten times.
We employ the GMRES solver that computes the residual.
The cost of multiplication of the $\mathbb{H}_{\epsilon}*x$ and the cost of multiplication of $\mathbb{A}*\mathbb{H}_{\epsilon}x$ is included in the fourth and fifth column in Table \ref{tab:input_2D}. The total cost of the GMRES with hierarchical matrices augmented by DNN compression is equal to {\bf compression cost of ${\mathbb H}_{\epsilon}$ with DNN} plus {\bf number of iterations} times the {\bf multiplication cost of ${\mathbb H}_{\epsilon}*x$ plus multiplication cost of ${\mathbb A}*{\mathbb H}_{\epsilon}x$}. The total cost is presented in the last column of Table \ref{tab:input_2D}.}

\inred{For comparison, we run the GMRES algorithm on the Galerkin method. The comparison is summarized in Table \ref{tab:input_2Da}. The number of iterations, the cost per iteration, and the total cost are presented there.
We can see that the cost of the stabilized solution is of the same order as the cost of the Galerkin solution. We compare here the costs of the correct solution obtained from the stable Petrov-Galerkin method with the cost of the incorrect solution from the unstable Galerkin method.}

\begin{figure}
    \centering

\includegraphics[width=0.9\textwidth]{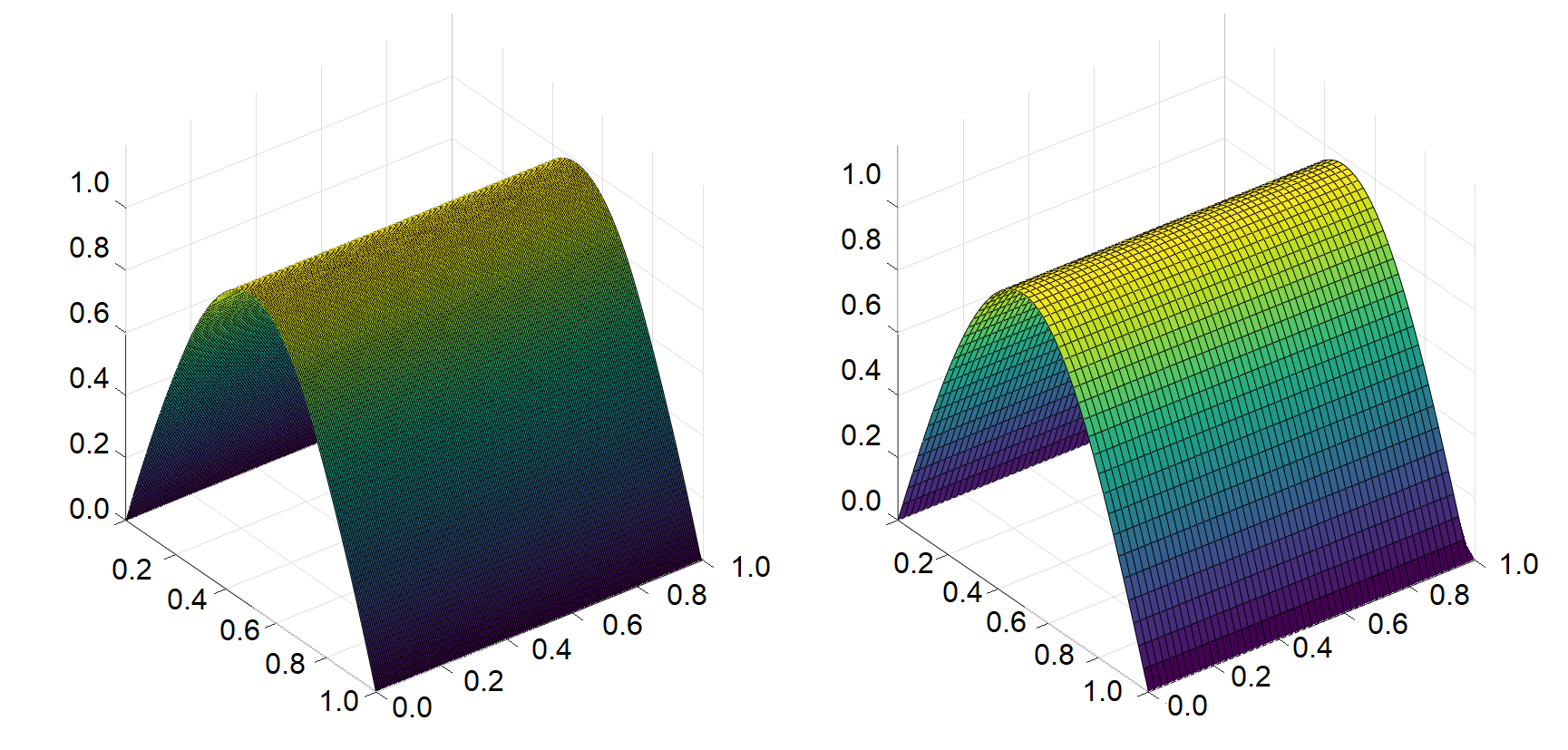}
    \includegraphics[width=0.9\textwidth]{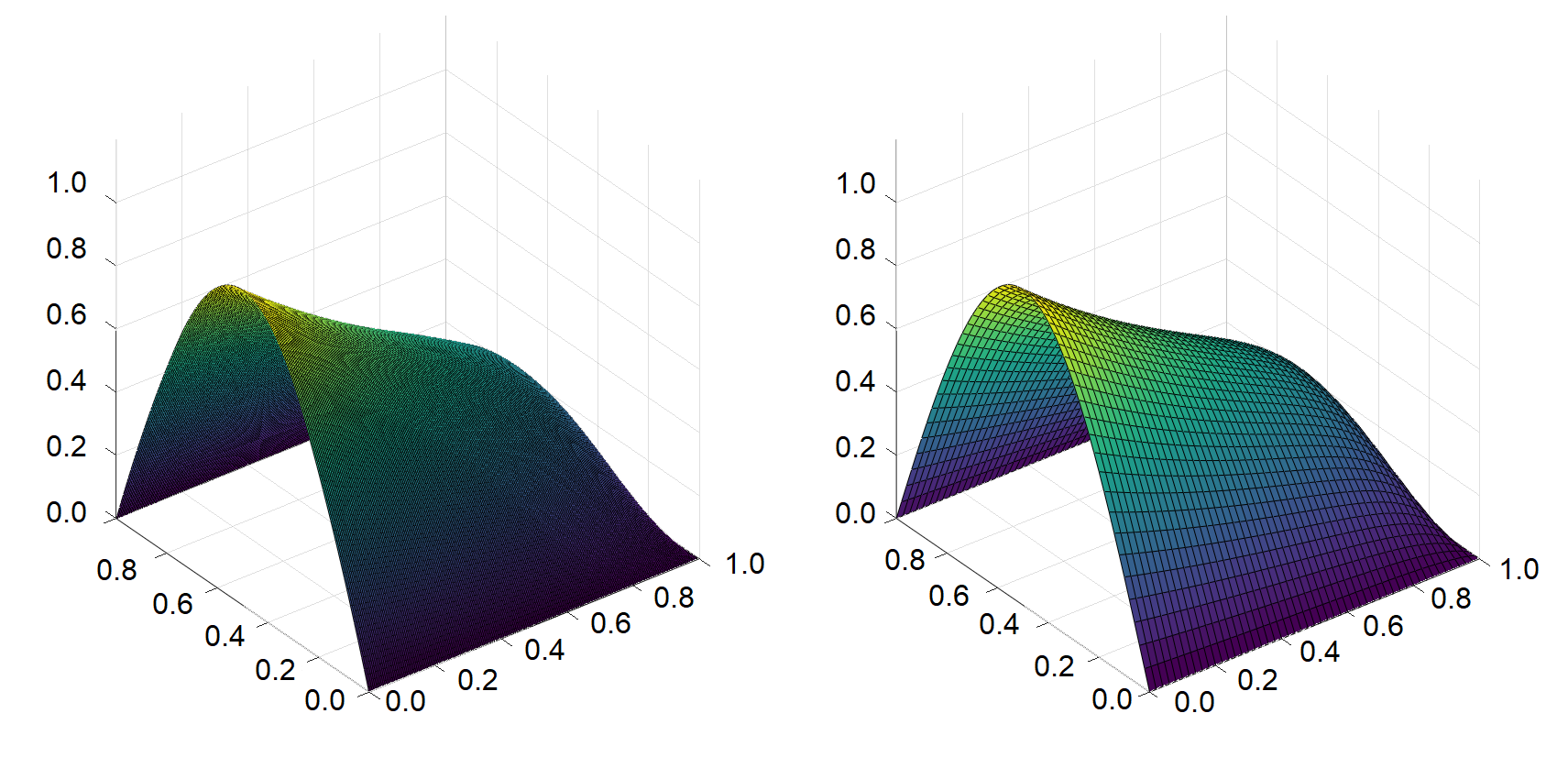}\\
        
    \caption{Case inflow data $g=\sin{\pi y}$. Comparison of exact solution (first column) and solutions obtained from GMRES solver using H-matrix (second column) for $\epsilon=0.000001$ (first row) and $\epsilon=0.1$ (second row). Iterative solver executed with accuracy $10^{-10}$.}
    
\label{fig:Hcomparison1}
\end{figure}
\begin{figure}
\center{
\includegraphics[scale=0.4]{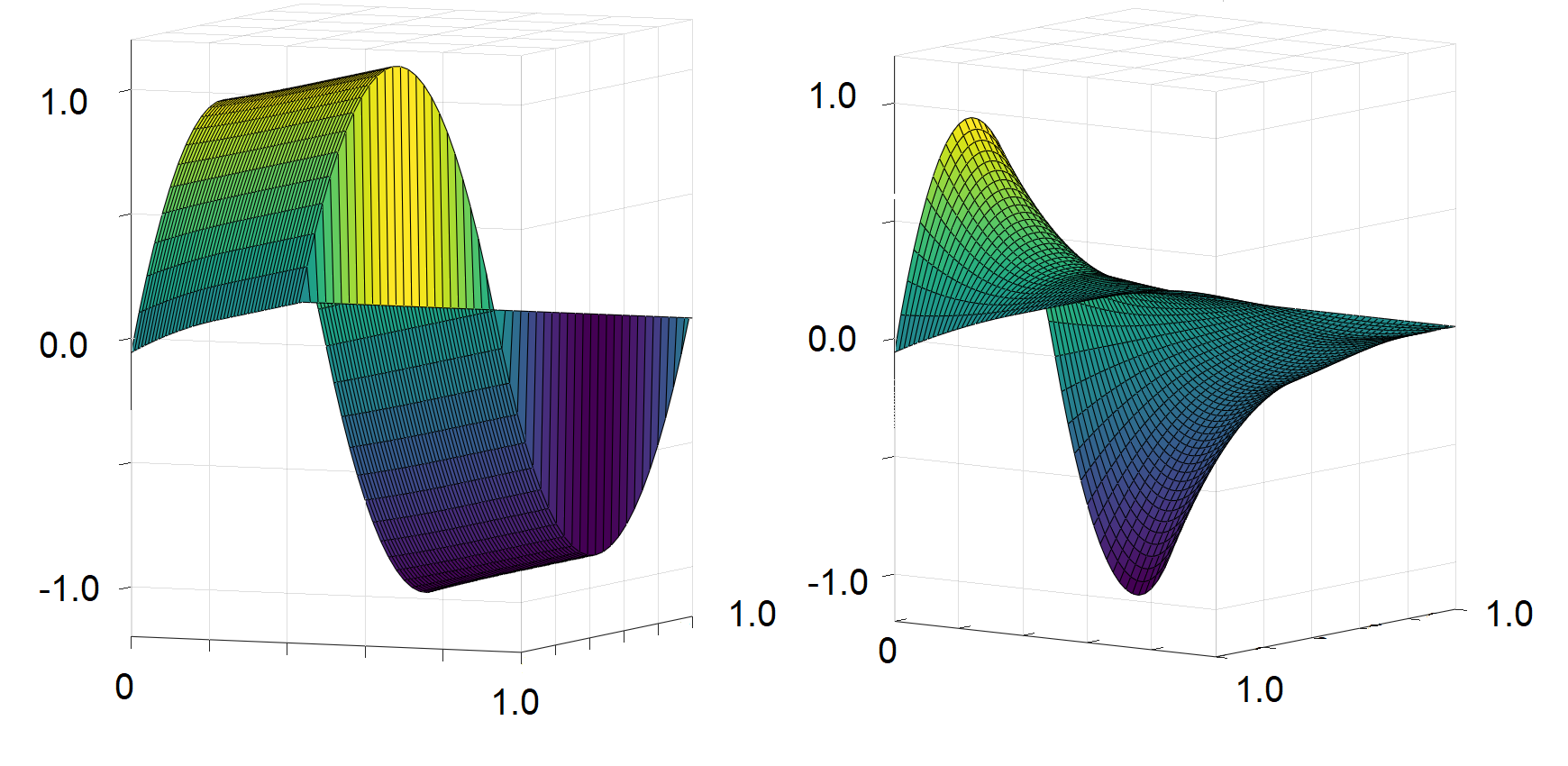}
}
\caption{Case inflow data $g=\sin(2\pi y)$. Solutions obtained from the Petrov-Galerkin formulation with  $\epsilon=0.000001$ and $\epsilon=0.1$.}
\label{fig:problem3}
\end{figure}

The numerical results are compared with the exact solutions in Figure \ref{fig:Hcomparison1} for inflow data $g=\sin(\pi y)$, and Figure~\ref{fig:problem3} for inflow data $g=\sin(2\pi y)$.




 \subsection{Hemholtz problem}
Given $\Omega=(0,1)^2\subset\mathbb R^2$ and $\kappa \in [1,10]$, we seek the solution of the Hemholtz problem

 \begin{equation}
 \left\{
 \begin{array}{rl}
 \Delta u + \kappa^2  u= f & \hbox{in }\Omega\\
 u=g & \hbox{over }\partial\Omega\,,
 \end{array}\right.
 \label{eq:Hemholtz}
 \end{equation}
with right-hand sides $f$ and $g$ such that the exact solution is
$u(x,y)=\sin(\kappa \pi x)\sin(\kappa \pi y)$.

We employ $20 \times 20$ finite elements mesh. The trial space is constructed from quadratic B-splines. The test space is obtained for trial, and quadraticfrom B-splines with $C^0$ separators (equivalent to Lagrange polynomials).
\inblue{The dependence of the coefficients of the optimal test functions on $\kappa$ for the Hemholtz problem has the affine structure described in Section~\ref{sec:affine}, thus we can \emph{offline}  construct the function}

\begin{equation}
\mathcal P\ni \kappa \rightarrow \mathbb W(\kappa),
\end{equation}

where $\mathbb W(\kappa)$ is the matrix of the coefficients of the optimal test functions.
We fix the trial and test spaces used for approximation of the solution and stabilization of the Petrov-Galerkin formulation. 
Next, we consider blocks of different size of matrix $\mathbb W(\kappa)$, and we train the SVD for these different blocks as function of $\kappa$, i.e., 
\begin{equation}
    \mathcal P\ni \kappa \rightarrow \operatorname{DNN}(\kappa)
    =\{\mathbb U_i(\kappa),\mathbb D_i(\kappa),\mathbb V_i(\kappa)\}_{i=1,...,N_B}\,.
\end{equation}

\begin{figure}
    \centering
\includegraphics[width=0.3\textwidth]{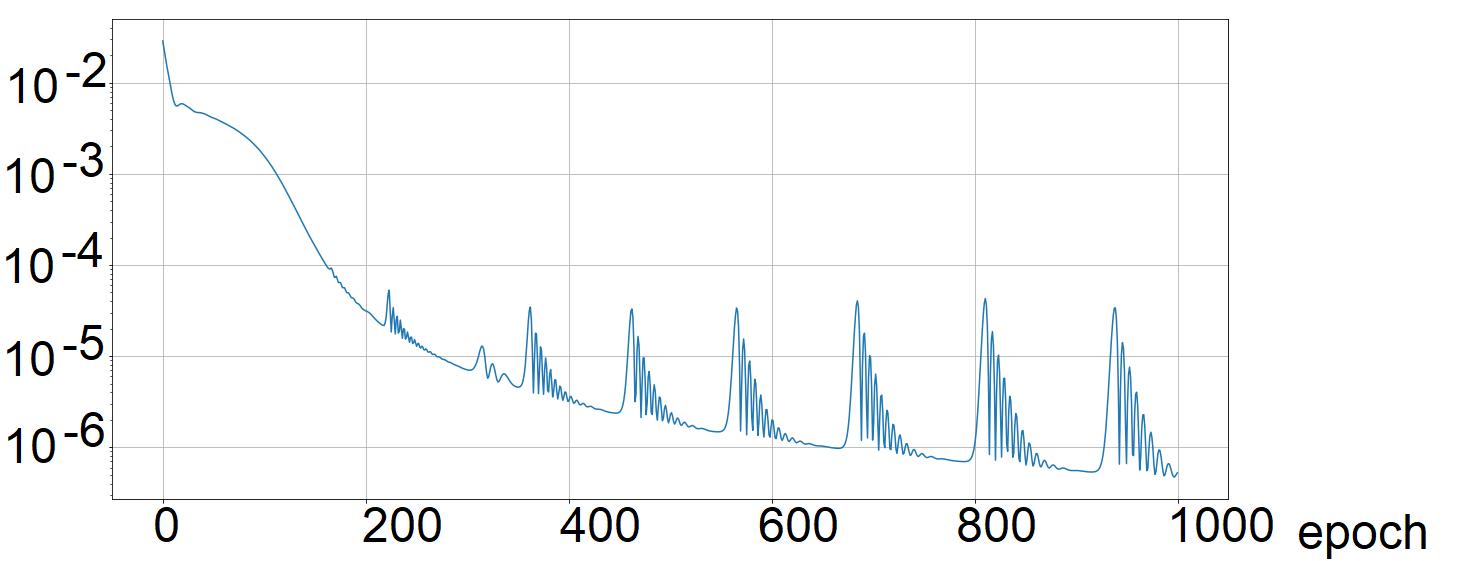}\includegraphics[width=0.3\textwidth]{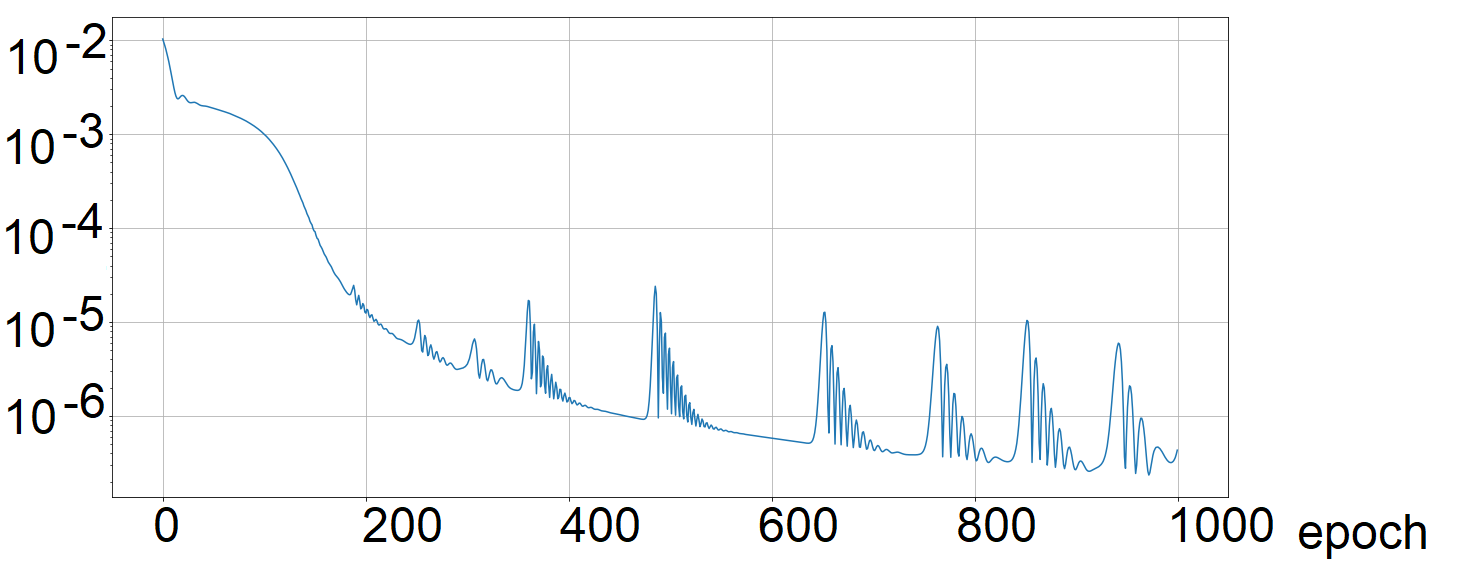}\includegraphics[width=0.3\textwidth]{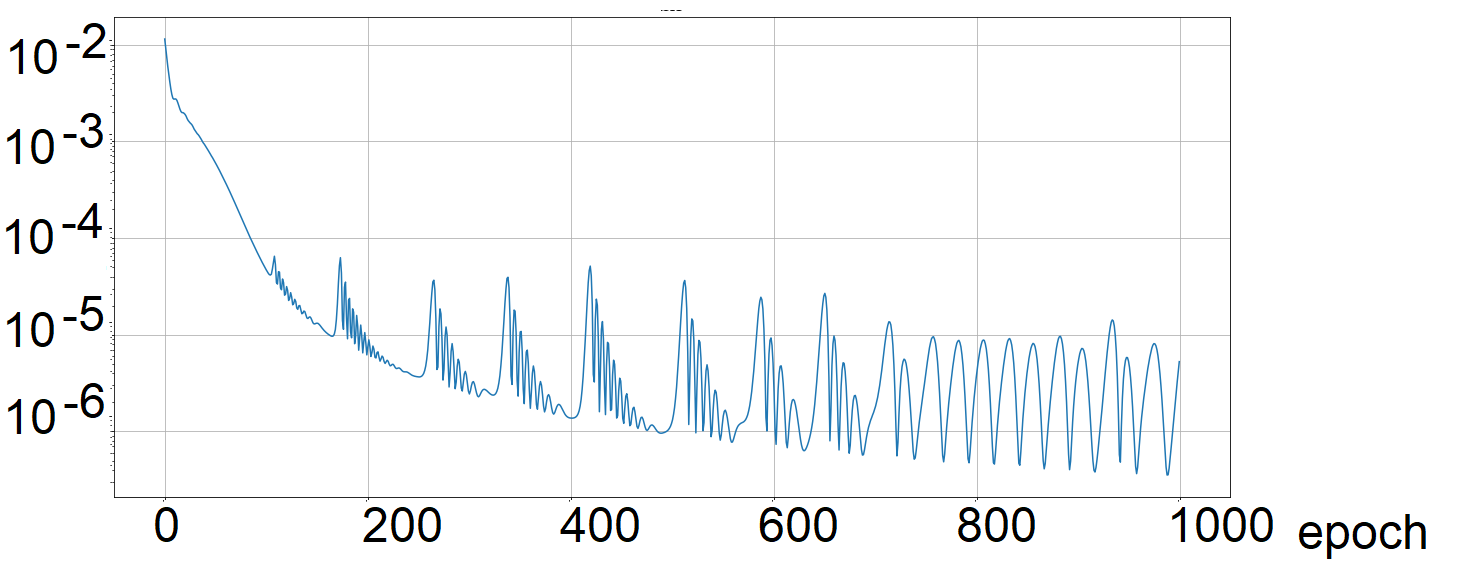}
    \caption{Hemholtz problem. Training of the neural network for the diagonal matrices $\mathbb U_{ij}(\kappa)$, $\mathbb D_{ij}(\kappa)$, and $\mathbb V_{ij}(\kappa)$ for one block.}
\label{fig:training1x1h}
\end{figure}

\inred{The convergence of the training procedure is presented in Figure \ref{fig:training1x1h}.
Knowing $\mathbb D_i(\kappa)$ a priori for a given $\kappa$ allows us to construct the structure of the hierarchical matrix, and we obtain the $\mathbb U_i(\kappa)$ and $\mathbb V_i(\kappa)$ from the neural networks}.

\begin{figure}
 \center{
 \includegraphics[scale=0.4]{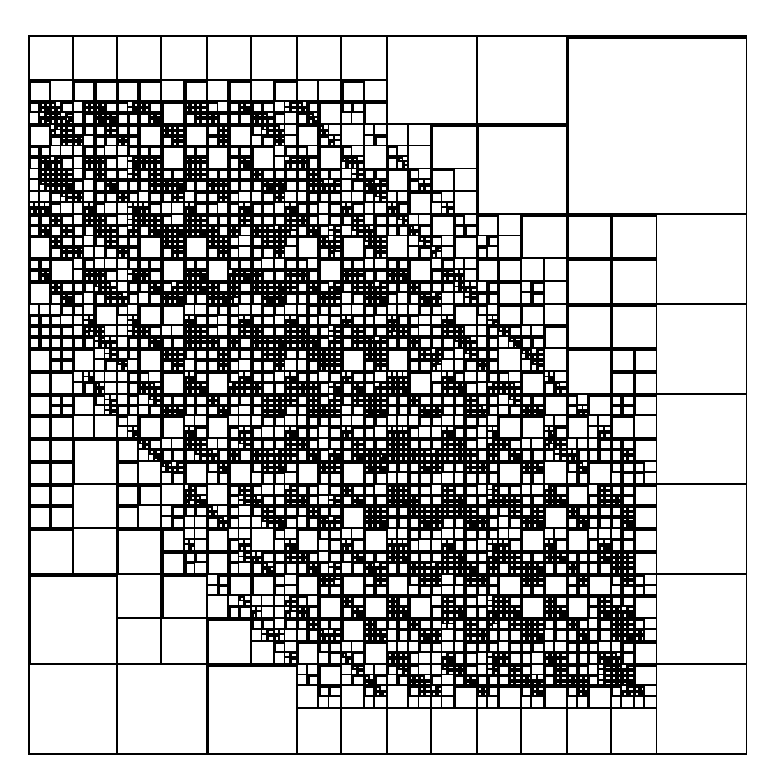}
 \includegraphics[scale=0.4] {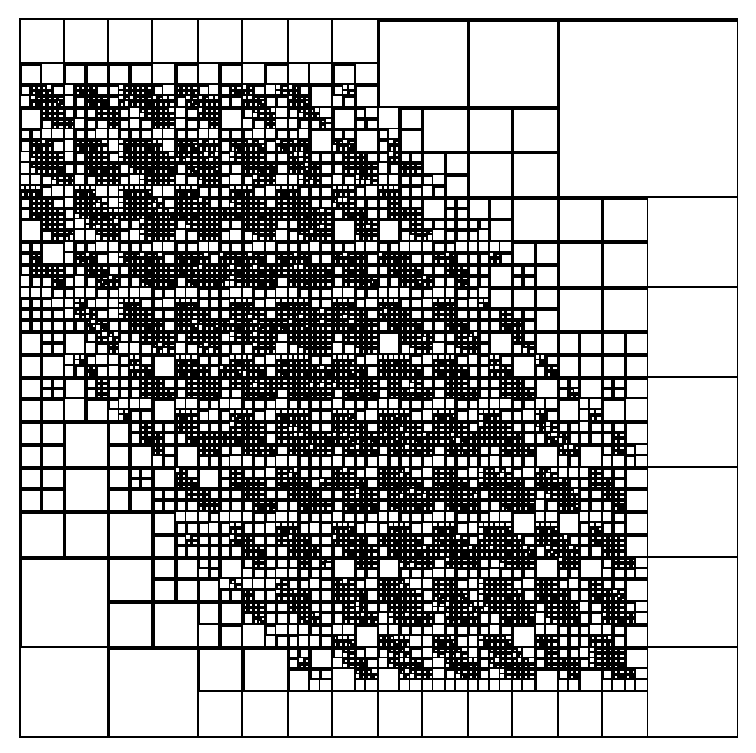}
 }
 \caption{The hierarchical matrices for $\kappa=1$ (left panel) and $\kappa=10$ (right panel). }
 \label{fig:hem_hmat}
\end{figure}

Figure \ref{fig:hem_hmat} depicts the exemplary resulting hierarchical matrices.

\inred{Table \ref{tab:input_2DH} summarizes the computational costs of our solver for two values of $\kappa=\{1,10\}$. The computational mesh was a tensor product of quadratic B-splines with 10 elements along $x$ and $y$ axes. The DNN allows obtaining the compression of the matrix of optimal test function's coefficients for free.
We employ the GMRES solver that computes the residual.
The cost of multiplication of the $\mathbb{H}_{\kappa}*x$ and the cost of multiplication of $\mathbb{A}*\mathbb{H}_{\kappa}x$ is included in the fourth and fifth column in Table \ref{tab:input_2DH}. The total cost of the GMRES with hierarchical matrices augmented by DNN compression is equal to the {\bf number of iterations} times the {\bf multiplication cost of ${\mathbb H}_{\kappa}*x$ plus multiplication cost of ${\mathbb A}*{\mathbb H}_{\kappa}x$}. The total cost is presented in the last column of Table \ref{tab:input_2DH}.}

\inred{We present in Table \ref{tab:input_2DHa} the cost of the GMRES algorithm executed on the Galerkin method. We present the number of iterations, the cost per iteration, and the total cost.
We can see that the cost of obtaining the stabilized solution is of the same order as the cost of the Galerkin method. We compare here the costs of the correct solution obtained from the Petrov-Galerkin method with the cost of the incorrect solution from the unstable Galerkin method.}

\begin{figure}
 \center{
 \includegraphics[scale=0.25]{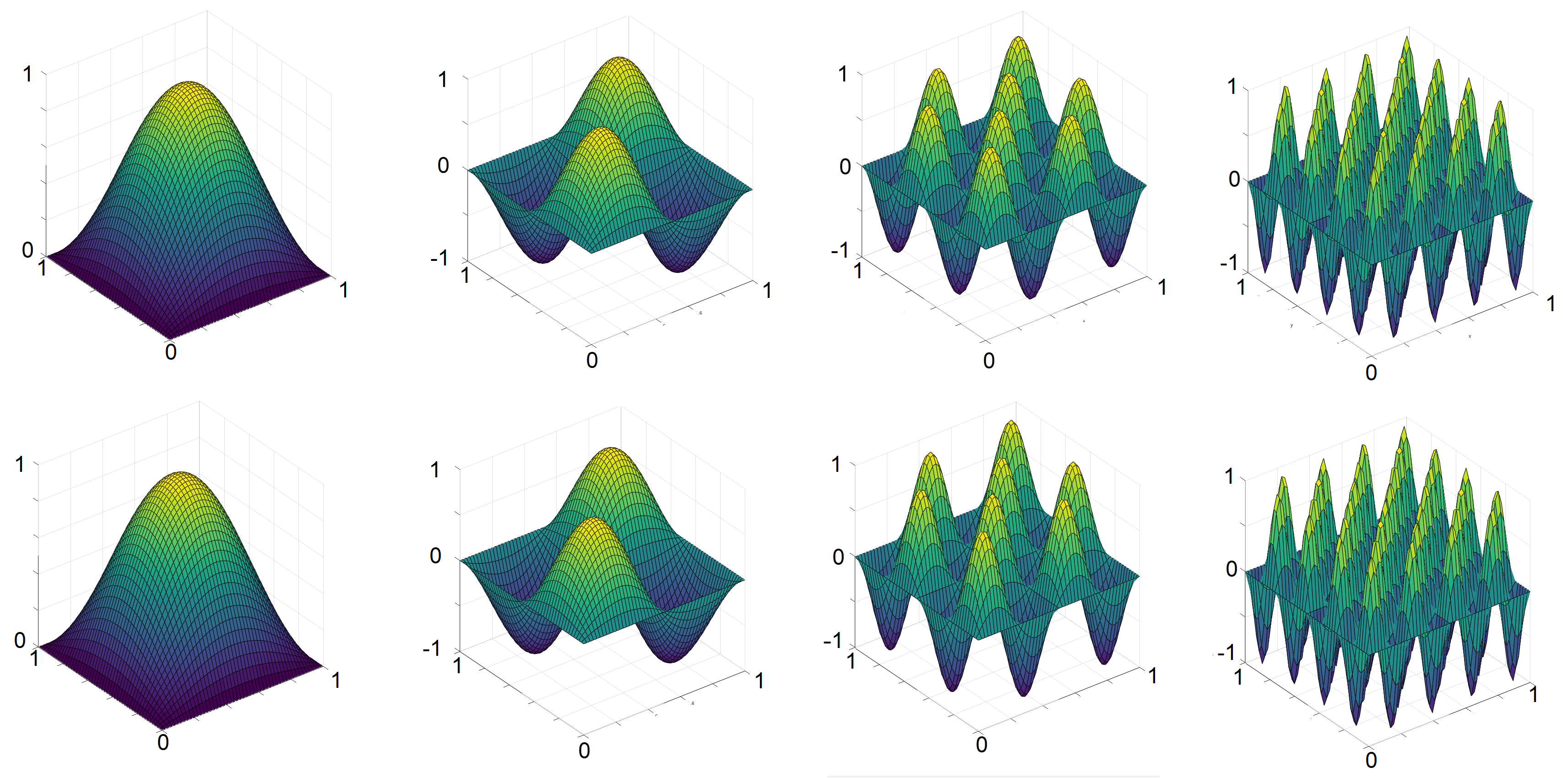}
 }
 \caption{The solutions obtained  for the Hemholtz problem for $\kappa=1,2,4,8$.
 The solution obtained from 
 the Petrov-Galerkin formulation with the optimal test functions provided by the neural network (first row).
 The exact solution (second row).}
 \label{fig:Hemholtz}
 \end{figure}

The comparison of the solution obtained with the Petrov-Galerkin formulation with the optimal test functions generated by DNN and the exact solution is presented in Figure \ref{fig:Hemholtz}.

\begin{table}
\begin{tabular}{ccccccc}
$\kappa$ & {\bf Compress } & {\bf Compress } & {\bf ${\mathbb W}_{\kappa}*x$} & {\bf ${\mathbb A}*{\mathbb W}_{\kappa}x$} & {\bf \# iter} & {\bf Total} \\
& {\bf ${\mathbb W}_{\kappa}$ flops} & {\bf ${\mathbb W}_{\kappa}$ flops} & {\bf flops} & {\bf flops} & {\bf GMRES } & {\bf flops}  \\
& {\bf with DNN} & {\bf without DNN} &  & & {\bf ${\mathbb H}$-matrix} & \\
 \hline
1 & 0 & 129,788 & 47,649 & 8,901 & 10 & 448,284 \\ 
10 & 0 & 129,788 & 47,649 & 8,877 & 31 & 1,752,306  \\
\hline
\vspace{5pt}
\end{tabular}
\caption{Computational costs of the stabilized Hemholtz solver using neural networks, hierarchical matrices and GMRES solver.}
\label{tab:input_2DH}
\end{table}

\begin{table}
\begin{tabular}{cccc}
$\epsilon$ &  {\bf \# iter GMRES Galerkin} & {\bf Flops per iteration} & {\bf Total flops}  \\
 \hline
1 &  10 & 14,444 & 144,440\\ 
10 & 31 & 14,444 & 447,764 \\
\hline
\vspace{5pt}
\end{tabular}
\caption{Computational costs of Galerkin method for the Hemholtz problem using GMRES solver.}
\label{tab:input_2DHa}
\end{table}

\section{Conclusions}
\inblue{We have employed the Petrov-Galerkin formulation with optimal test functions for the stabilization of difficult problems. We have focused on advection-dominated diffusion and Helmholtz problems.
During the \emph{offline} phase, we explicitly compute the matrix of coefficients of optimal test functions for any PDE-parameter. We have also trained neural networks to compute for each PDE-parameter the bottleneck of hierarchical matrix compression.
During the \emph{online} phase, we rapidly compute the matrix compression using the neural networks, and we perform the GMRES iterative solver on the reduced Petrov-Galerkin linear system, where vector-matrix multiplications are done in a quasi-linear computational cost, due to the hierarchical structure of the low-rank decomposition used. 
Thus, we obtain the \emph{online} stabilization practically for free.}

\section*{Acknowledgments}
The European   Union's Horizon 2020 Research and Innovation Program of the Marie Sk\l{}odowska-Curie grant agreement No. 777778, MATHROCKs.
Research project partly supported by the program ``Excellence initiative – research university" for the University of Science and Technology.

\appendix
\section{Algorithms}
\subsection{Recursive hierarchical compression of the matrix}
\begin{algorithm}[H]
\begin{algorithmic}[h]
\REQUIRE {$t_{\min},t_{\max},s_{\min},s_{\max} \in {\mathbb {N}}$ (row and column index ranges),
\newline
$1 \le t_{\min} \le t_{\max} \le n,1 \le s_{\min} \le s_{\max} \le m $ where $n \times m$ is the size of the matrix to be compressed} 
\IF {$\operatorname{Admissible}(t_{\min},t_{\max},s_{\min},s_{\max},r,\delta)$}
\STATE $v=\operatorname{CompressMatrix}(t_{\min},t_{\max},s_{\min},s_{\max},r)$
\ELSE
\STATE // Create a new node with four sons corresponding to four //quarters of the matrix 
\STATE create new node $v$
\STATE AppendChild(v,CreateTree$(t_{\min},t_{\operatorname{newmax}},s_{\min}, s_{\operatorname{newmax}}))$
\STATE AppendChild(v,CreateTree$(t_{\min},t_{\operatorname{newmax}},s_{\operatorname{newmax}}+1, s_{\max}))$
\STATE AppendChild(v,CreateTree$(t_{\operatorname{newmax}}+1,t_{\max},s_{\min}, s_{\operatorname{newmax}}))$
\STATE AppendChild(v,CreateTree$(t_{\operatorname{newmax}}+1,t_{\max},s_{\operatorname{newmax}}+1, s_{\max}))$
\ENDIF
\STATE RETURN $v$
\end{algorithmic}
\caption{Recursive hierarchical compression of the matrix: \textbf{CreateTree($r,\delta$)} where $r$ is the rank used for the compression, and $\delta$ is the threshold for the 
singular values.}
\label{Alg1}
\end{algorithm}
\subsection{Checking of the admissibility condition}
\begin{algorithm}[H]
\begin{algorithmic}[h]
\REQUIRE $t_{\min},t_{\max},s_{\min},s_{\max}$ - range of indexes of block,
$\delta$ compression threshold, $r$ maximum rank
\IF {block ($t_{\min},t_{\max},s_{\min},s_{\max}$) consist of zeros}
\STATE \textrm{\bf return true};
\ENDIF
\STATE {$[\mathbb U,\mathbb D,\mathbb V]\gets \operatorname{truncatedSVD}(t_{\min},t_{\max},s_{\min},s_{\max},r+1)$; $\sigma \gets \operatorname{diag}(\mathbb D)$;}
\IF {$\sigma(r+1)<\delta$}
\STATE {\bf return true};
\ENDIF
\STATE {\bf return false};
\end{algorithmic}
\caption{Checking of the admissibility condition: $\operatorname{result} = \operatorname{Admissible}(t_{\min},t_{\max},s_{\min},s_{\max},r,\delta)$}
\label{Alg2}
\end{algorithm}
\subsection{Matrix vector multiplication}
\begin{algorithm}[H]
\caption{\revision{Matrix vector multiplication: $Y = {\bf matrix\_vector\_mult}(v,X)$}}
\begin{algorithmic}[h]
\REQUIRE $\textrm{node }v$  \ingreen{ representing} \emph{compressed matrix} $\mathbb H(v) \in {\cal M}^{m\times n}$,  $X \in {\cal M}^{n \times c}$ \emph{vectors to multiply}
\IF {$v.nr\_sons==0$}
\IF {$v.rank==0$}
\STATE $\textrm{return }zeros(size(A).rows)$
\ENDIF
\STATE $\textrm{return }v.U*(v.V*X)$
\ENDIF
\STATE $rows = size(X).rows$
\STATE $X_1=X(1:\frac{rows}{2},*)$
\STATE $X_2=X(\frac{rows}{2}+1:size(A).rows,*)$
\STATE $C_2=v.son(1).U; C_1=v.son(1).V$
\STATE $D_2=v.son(2).U; D_1=v.son(2).V$
\STATE $E_2=v.son(3).U; E_1=v.son(3).V$
\STATE $F_2=v.son(4).U; F_1=v.son(4).V$
\STATE $\textrm{return } \begin{bmatrix} C_2*(C_1*X_1)+D_2*(D_1*X_2) \\  E_2*(E_1*X_1)+F_2*(F_1*X_2) \end{bmatrix}$
\end{algorithmic}
\end{algorithm}
\subsection{rSVD compression of a block}
\begin{algorithm}[H]
\begin{algorithmic}[h]
\REQUIRE $t_{\min},t_{\max},s_{\min},s_{\max}$ - range of indexes of block,
$\delta$ \emph{compression threshold}, $r$ \emph{maximum rank}
\IF {block ($t_{\min},t_{\max},s_{\min},s_{\max}$) consist of zeros}
\STATE $\textrm{\bf create new node } v; v.rank \gets 0; v.size \gets size(t_{\min},t_{\max},s_{\min},s_{\max}); \textrm{\bf return }v;$
\ENDIF
\STATE $[\mathbb U,\mathbb D,\mathbb V]\gets \operatorname{truncatedSVD}(t_{\min},t_{\max},s_{\min},s_{\max},r)$; $\sigma \gets \operatorname{diag}(\mathbb D)$;
\STATE $rank \gets rank(\mathbb D)$
\STATE $\textrm{\bf create new node }v;$ $v.rank \gets rank;$
\STATE $v.singularvalues \gets \sigma(1:rank);$
\STATE $v.U \gets U(*,1:rank);$
\STATE $v.V \gets D(1:rank,1:rank)*V(1:rank,*);$
\STATE $v.sons \gets\emptyset;$ $v.size \gets size(t_{\min},t_{\max},s_{\min},s_{\max});$
\STATE $\textrm{\bf return }v;$
\end{algorithmic}
\caption{rSVD compression of a block: $node = {\bf CompressMatrix}(t_{\min},t_{\max},s_{\min},s_{\max},r)$}
\label{Alg3}
\end{algorithm}
\subsection{Pseudo-code of the GMRES algorithm}
\begin{algorithm}[H]
\begin{algorithmic}[h]
\REQUIRE $A$ matrix, $b$ right-hand-side vector, $x_0$ starting point
\STATE Compute $r_0=b-Ax_0$
\STATE Compute $v_1=\frac{r_0}{\|r_0\|}$
\STATE {\bf for}
$j = 1, 2, . . . , k$
\STATE 
Compute $h_{i,j} = \left(Av_j , v_i\right)$
 for $i = 1, 2, . . . , j$
\STATE Compute $\hat{v}_{j+1}=Av_j-\sum_{i=1,...,j}h_{i,j}v_i$
\STATE Compute $h_{j+1,j} = \|\hat{v}_{j+1}\|_2$
\STATE Compute $v_{j+1}=\hat{v}_{j+1}/h_{j+1,j}$
\STATE {\bf end for}
\STATE Form solution $x_k=x_0+V_ky_k$ where 
$V_k=[v_1 ... v_k]$,
and $y_k$ minimizes
$J(y)=\|\beta e_1-\hat{H}_ky\|$ where
$\hat{H}=\begin{bmatrix}h_{1,1} & h_{1,2} \cdots h_{1,k}\\
h_{2,1} & h_{2,2} \cdots h_{2,k}\\
0 &  \ddots & \ddots & \vdots \\
\vdots & \ddots & h_{k,k-1} & h_{k,k} \\
0 & \cdots & 0 & h_{k+1,k}
\end{bmatrix}$
\label{Alg4}
\end{algorithmic}
\caption{Pseudo-code of the GMRES algorithm}
\end{algorithm}
\subsection{Recursive hierarchical compression of the matrix augmented by neural network}
\begin{algorithm}[H]
\begin{algorithmic}[h]
\REQUIRE {$t_{\min},t_{\max},s_{\min},s_{\max}, \in {\mathbb {N}}$ (row and column index ranges), $r$ rank of the blocks, $\delta$ accuracy of compression, $\mu$ PDE parameter
\newline
$1 \le t_{\min} \le t_{\max} \le n,1 \le s_{\min} \le s_{\max} \le m $ where $n \times m$ is the size of the matrix to be compressed} 
\STATE i = block index for $(t_{\min},t_{\max},s_{\min},s_{\max})$
 \IF { $\mathbb D_i(\mu)[r+1]<\delta$  \emph{(asking NN for block singularvalues}}
 \IF {block ($t_{\min},t_{\max},s_{\min},s_{\max}$) consist of zeros}
\STATE $\textrm{\bf create new node } v; v.rank \gets 0; v.size \gets size(t_{\min},t_{\max},s_{\min},s_{\max},s,t); \textrm{\bf return }v;$
\ENDIF
\STATE $[\mathbb U,\mathbb D,\mathbb V]\gets \operatorname{truncatedSVD}(t_{\min},t_{\max},s_{\min},s_{\max},r); \sigma \gets \operatorname{diag}(\mathbb D);$
\STATE $rank \gets rank(\mathbb D)$
\STATE $\textrm{\bf create new node }v;$ $v.rank \gets rank;$
\STATE $v.singularvalues \gets \sigma(1:rank);$
\STATE $v.U \gets \mathbb U(*,1:rank);$
\STATE $v.V \gets \mathbb D(1:rank,1:rank)*\mathbb V(1:rank,*);$
\STATE $v.sons \gets\emptyset;$ $v.size \gets size(t_{\min},t_{\max},s_{\min},s_{\max});$
\STATE $\textrm{\bf return }v;$
\ELSE
\STATE // Create a new node with four sons corresponding to four //quarters of the matrix 
\STATE create new node v
\STATE AppendChild(v,CreateTreeNN$(t_{\min},t_{\operatorname{newmax}},s_{\min}, s_{\operatorname{newmax}},r,\delta,\mu))$
\STATE AppendChild(v,CreateTreeNN$(t_{\min},t_{\operatorname{newmax}},s_{\operatorname{newmax}}+1, s_{\max},r,\delta,\mu))$
\STATE AppendChild(v,CreateTreeNN$(t_{\operatorname{newmax}}+1,t_{\max},s_{\min}, s_{\operatorname{newmax}},r,\delta,\mu))$
\STATE AppendChild(v,CreateTreeNN$(t_{\operatorname{newmax}}+1,t_{\max},s_{\operatorname{newmax}}+1, s_{\max},r,\delta,\mu))$
\ENDIF
\STATE RETURN $v$
\end{algorithmic}
\caption{Recursive hierarchical compression of the matrix augmented by neural network: \textbf{CreateTreeNN(1,rowsof($A$),1,columnsof($A$),$r,\delta,\mu$)} where $r$ is the rank used for the compression, and $\delta$ is the threshold for the $r$  singular values, $\mu$ is the PDE parameter.}
\label{Alg2}
\end{algorithm}


\end{document}